\documentclass[preprint,11pt,number]{elsarticle}
\usepackage{setspace}
\usepackage[
    a4paper,
    left=2cm,
    right=2cm,
    top=2.5cm,
    bottom=2.5cm
]{geometry}

\usepackage{amssymb}
\usepackage{amsmath}
\usepackage{amsmath}
\usepackage{amsthm}
\usepackage{amssymb}
\usepackage{mathtools}
\usepackage{graphicx}
\usepackage{subcaption}
\usepackage[colorlinks=true, allcolors=blue]{hyperref}
\usepackage{tikz}
\usepackage{booktabs}
\usepackage{multirow}
\usepackage{algorithm}
\usepackage{algpseudocode}
\usetikzlibrary{positioning,arrows.meta, fit, calc}
\usepackage{float}
\usepackage{placeins}
\usepackage{makecell}
\usetikzlibrary{decorations.pathmorphing}
\usepackage{booktabs}
\usepackage{siunitx}
\usepackage{array}
\usepackage[dvipsnames]{xcolor}
\usepackage{enumitem}
\usepackage{rotating}
\usepackage{booktabs} 
\usepackage{pifont}   
\usepackage{pgfplots}
\usepackage{threeparttable}
\usepgfplotslibrary{groupplots}
\newcommand{\cmark}{\ding{51}} 
\newcommand{\xmark}{\ding{55}} 
\newcolumntype{C}{>{\centering\arraybackslash}X}
\newcolumntype{R}{>{\raggedleft\arraybackslash}X}

\newtheorem{lemma}{Lemma}

\journal{European Journal of Operational Research}
\begin{document}

\begin{frontmatter}



\title{Integrated Surgical Scheduling with Weekend Bed Occupancy Management: A Column Generation Approach} 


\author[p]{Sara Cambiaghi\corref{cor1}}
\cortext[cor1]{Corresponding author}
\ead{sara.cambiaghi01@universitadipavia.it}
\author[l]{Saverio Basso}
\author[p]{Davide Duma}
\author[l]{Matteo Salani}

\affiliation[l]{organization={SUPSI, Dalle Molle Institute for Artificial Intelligence (IDSIA)},
addressline={Via la Santa 1}, 
city={Lugano},
postcode={6900},
country={Switzerland}}

\affiliation[p]{organization={Department of Mathematics “F. Casorati”, University of Pavia},
addressline={Via Ferrata, 5}, 
city={Pavia},
postcode={27100},
country={Italy}}

\begin{abstract}
This paper addresses the integrated Master Surgical Scheduling Problem and Surgical Case Assignment Problem under weekend bed capacity constraints. We propose a joint optimization framework that simultaneously determines the assignment of specialties and surgeons to Operating Room (OR) blocks, the selection and assignment of patients to these blocks, and the sequencing of surgeries, while accounting for reduced bed availability during weekends. The model is formulated as a multi-criteria Mixed-Integer Linear Program (MILP) that prioritizes the total priority-weighted amount of scheduled patients, minimizes excess weekend bed occupancy, and maximizes OR utilization. To solve large instances, we develop a Column Generation (CG) algorithm in which the pricing subproblem is formulated as a Resource-Constrained Shortest Path Problem. The framework accommodates block, open, and modified block scheduling policies. Computational experiments show that the CG approach outperforms the direct MILP solution in terms of solution quality and computational time, particularly for larger instances. We further analyze the trade-offs between weekend bed availability, OR utilization, and total patient priority. Our results provide managerial insights into how hospitals can balance patient access, operational efficiency, and staff workload when resources are limited.
\end{abstract}


\begin{highlights}
\item A multi-criteria optimization model integrating the master surgical scheduling problem, surgical case assignment, and case sequencing.
\item A column generation framework with the subproblem formulated as a resource-constrained shortest path problem.
\item Explicit modeling of weekend bed capacity constraints to account for reduced staffing levels.
\item A comparative analysis of block, open, and modified block scheduling policies.
\item A sensitivity analysis of resource availability providing managerial insights.
\end{highlights}

\end{frontmatter}


\section{Introduction}
\label{section:introduction}

Managing Operating Rooms (ORs) is a critical task in hospitals, both from an efficiency perspective and in terms of patient satisfaction and safety. From a financial standpoint, the operating theater generates more than 40\% of a hospital’s revenue while accounting for nearly 30\% of its total expenditures \cite{AlAmin}. From the patient perspective, population aging is increasing the demand for health services, leading to longer waiting times and placing growing pressure on hospitals \cite{Aringhieri}.

OR planning is inherently complex because it involves multiple stakeholders, including patients, surgeons, nurses, and other clinical staff, as well as several interdependent resources, such as ORs and inpatient beds. Consequently, most existing studies primarily focus on the allocation of OR capacity, whereas downstream resources such as inpatient beds receive comparatively less attention. However, as highlighted by Landa et al.\ \cite{Landa}, inpatient beds are also a scarce resource whose capacity cannot realistically be adjusted in the short term. Consequently, bed availability represents a hard constraint and often becomes the main bottleneck in OR scheduling \cite{Mihalj2022}.

This issue is particularly relevant during weekends, when hospital staffing levels are typically reduced. Although beds remain physically available, admitting patients up to the nominal bed capacity is generally undesirable, as the reduced availability of medical and nursing staff may limit the hospital's ability to provide an adequate level of care \cite{Bell}. For this reason, hospitals usually operate with lower target occupancy levels during weekends, reflecting the fact that the effective capacity of a hospital depends not only on the number of available beds but also on the resources required to provide appropriate treatment and observation \cite{Bechir}.

An operating theater (OT) consists of a set of ORs within the same hospital unit, serving one or more surgical specialties or hospital wards. Since these specialties share the same OR capacity, their scheduling decisions are inherently interdependent. Consequently, surgical planning involves two closely related decision problems.

The first is the \textit{Master Surgical Scheduling Problem (MSSP)}, which determines the Master Surgical Schedule (MSS), namely the allocation of OR blocks to surgical specialties over a given planning horizon. An OR block represents the availability of a specific OR during a given time period, typically a day or part of a day. The resulting MSS defines the distribution of OR capacity among specialties and provides the framework for subsequent patient scheduling decisions.

According to Zhu et al.~\cite{Zhu}, hospitals typically adopt one of three scheduling policies. In \textit{block scheduling}, OR blocks are assigned in advance to specific specialties, simplifying the planning process. In \textit{open scheduling}, OR blocks are shared across specialties, providing greater scheduling flexibility. Finally, \textit{modified block scheduling} combines dedicated and shared OR blocks, balancing planning simplicity with scheduling flexibility. 
According to Augusto et al.~\cite{Augusto}, this policy is commonly implemented either by reserving only a subset of OR blocks for specific specialties while keeping the remaining blocks shared, or by allowing unused capacity within specialty-dedicated blocks to be reassigned to other specialties once a predefined utilization threshold is exceeded.

The second problem is the \textit{Surgical Case Assignment Problem (SCAP)}, which consists of assigning elective patients to OR blocks and determining the sequence of surgeries within each block while accounting for OR capacity, bed availability, and surgeon availability.

In this work, we consider a set of surgical specialties and surgeons, together with a set of patients characterized by their specialty, surgery duration, expected Length of Stay (LoS), assigned surgeon, and a priority score. Since not all patients can be scheduled within the planning horizon, a subset must be selected according to clinical priority and assigned to appropriate OR blocks. 
In addition, we define, for each specialty, a target level of weekend bed occupancy that is lower than the corresponding weekday capacity, reflecting the reduced staffing levels typically observed during weekends.
This modeling choice is motivated by the well-known \emph{weekend effect}, namely the poorer outcomes observed among patients admitted or treated during weekends across a wide range of medical settings, diagnoses, and countries. Although its underlying causes remain debated, the two most widely accepted explanations are reduced weekend staffing and the higher severity of patients admitted during weekends \cite{Honeyford}. 

In this work, surgeries are scheduled during weekdays while explicitly limiting weekend bed occupancy. Exceeding the target occupancy may place additional pressure on the reduced weekend staff. Therefore, deviations from the target are penalized, allowing limited flexibility while discouraging excessive weekend occupancy. 
We also investigate the trade-off between increasing weekend bed capacity, which requires additional staffing resources, and the resulting increase in the number of surgeries that can be performed during the week.

Although the MSSP, the SCAP, and downstream bed management have all been extensively studied in the literature, these problems are often addressed separately. As a result, the interactions between OR allocation, patient assignment, and downstream bed occupancy are only partially captured. Moreover, despite the operational relevance of reduced weekend staffing, weekend bed occupancy has received limited attention in integrated OR scheduling models.
This limitation motivates the integrated framework proposed in this paper

The proposed approach is based on a Column Generation (CG) algorithm and a multi-objective optimization model that (i) maximizes the total priority-weighted number of scheduled patients, (ii) minimizes excess weekend bed occupancy, and (iii) maximizes OR utilization. The resulting formulation promotes patient prioritization, efficient use of hospital resources, and improved coordination between surgical activity and downstream bed availability.

Finally, we compare block, open, and modified block scheduling policies under two representative operational settings, characterized by short and long patient LoS, respectively. We further evaluate their performance under different levels of resource availability, providing managerial insights into the trade-offs between patient access, OR utilization, and weekend bed capacity.

\section{Literature Review}
\label{section:literature_review}

In recent decades, OR scheduling has attracted considerable attention, resulting in numerous studies addressing different aspects of the problem and proposing a wide variety of models, often motivated by real hospital case studies. Several literature reviews provide comprehensive overviews of the field \cite{AlAmin, Cardoen, Guerriero, Gur, Hof, Nasirian, Rahimi, Razali, Samudra, VanRiet, Zhu}.

Razali et al.~\cite{Razali} distinguish between cyclic and non-cyclic Master Surgical Schedules (MSSs). Their review shows that cyclic schedules, which repeat after a fixed planning horizon (typically one week), are the most widely adopted because of their simplicity. However, their periodic structure may reduce flexibility compared with schedules that are re-optimized regularly to accommodate changes in patient demand and resource availability \cite{Oliveira2022}.

Taken together, these reviews show that OR scheduling models typically aim to maximize patient priority and OR utilization while minimizing costs, overtime, and inefficient resource allocation. Surgery duration, patient LoS, emergency arrivals, and resource availability are the main sources of uncertainty, which are commonly addressed through stochastic programming, robust optimization, simulation, and fuzzy optimization \cite{Nasirian}. During schedule execution, online optimization techniques can be used to adapt the schedule to unforeseen events \cite{Duma}.

\subsection{Integrated Master Surgical Scheduling, Surgical Case Assignment, and Bed Management}

This section focuses on studies that jointly address the MSSP and SCAP while also considering downstream bed management. To the best of our knowledge, only two studies explicitly integrate all three problems simultaneously.

Aringhieri et al.~\cite{Aringhieri} propose a multi-neighborhood local-search matheuristic for the joint MSSP and SCAP. Their objectives are to maximize patient priority while balancing bed occupancy across specialties through a bed-leveling strategy.

Siqueira et al.~\cite{Siqueira} develop an optimization model covering the entire patient pathway, from admission to discharge. Their approach jointly determines a periodic MSS and a specialty-specific bed allocation policy, and compares block and open scheduling policies, showing that block scheduling requires substantially more capacity.

Several other studies address subsets of these decisions (Table~\ref{tab:literature_comparison}). Testi et al.~\cite{Testi,Testi2} propose hierarchical approaches combining MSS design with patient assignment and admission planning, including a week-surgery setting in which all patients must be discharged before the weekend. Makboul et al.~\cite{Makboul2,Makboul3} investigate integrated MSSP-SCAP formulations using either multi-objective optimization or decomposition strategies under downstream bed-capacity constraints. Agnetis et al.~\cite{Agnetis} decompose the problem into MSS design and independent knapsack problems for patient assignment, substantially reducing computational time.

Other studies primarily focus on balancing downstream bed demand while constructing the MSS. Kianfar et al.~\cite{Kianfar} incorporate patient LoS, surgeon preferences, and bed occupancy into the MSS design. Similarly, Deklerck et al.~\cite{Deklerck} and Calegari et al.~\cite{Calegari} optimize cyclic MSSs to smooth ward bed occupancy, using multi-objective optimization and simulation-based analyses to evaluate solution quality under uncertainty.

\begin{sidewaystable}[p]
\centering
\caption{Comparison of the most closely related studies according to the main modeling assumptions and solution approaches (in chronological order). \textit{MSSP} indicates whether the Master Surgical Schedule problem is considered. \textit{SCAP} indicates whether the Surgical Case Assignment Problem is considered. \textit{Joint} indicates whether the MSS and SCAP are solved simultaneously. \textit{Beds} indicates whether bed availability is explicitly considered. \textit{Surgeons} indicates whether surgeon-related constraints are considered and how. \textit{Turnover} indicates whether turnover times are explicitly modeled and how. \textit{Sequencing} indicates whether the sequencing of surgical procedures within operating rooms is considered. \textit{Policy} indicates the operating room scheduling policy considered (block, open, or modified). \textit{Solution approach} indicates the optimization and solution methods adopted.}
\label{tab:literature_comparison}
\newcolumntype{C}[1]{>{\centering\arraybackslash}p{#1}}
\renewcommand{\arraystretch}{1.2}
\small
\resizebox{\textheight}{!}{%
\begin{tabular}{@{}lccccC{3cm}C{3cm}ccC{3cm}@{}}
\toprule
\textbf{Study} &
\textbf{MSSP} &
\textbf{SCAP} &
\textbf{Joint} &
\textbf{Beds} &
\textbf{Surgeons} &
\textbf{Turnover} &
\textbf{Sequencing} &
\textbf{Policy} &
\textbf{Solution approach} \\
\midrule

Testi et al.~(2007) \cite{Testi} 
& \cmark & \cmark & \xmark & \cmark
& Surgical team availability
& \xmark
& \cmark
& Block
& ILP + discrete-event simulation \\

Testi \& Tànfani (2009) \cite{Testi2}
& \cmark & \cmark & \cmark & \cmark
& Surgical team availability
& \xmark
& \xmark
& Block
& ILP \\

Agnetis et al.~(2014) \cite{Agnetis} 
& \cmark & \cmark & \xmark & \xmark
& Surgical team availability
& Included in procedure durations
& \cmark
& Block
& Minimum-cost flow + multiple knapsack \\

Siquiera et al.~(2018) \cite{Siqueira} 
& \cmark & \cmark & \cmark & \cmark
& \xmark
& Clean-up and set-up time
& \cmark
& Block/Open
& ILP \\

Makboul et al.~(2020) \cite{Makboul2} 
& \cmark & \cmark & \cmark & \xmark
& OR block assignment
& \xmark
& \cmark
& Block
& MILP + $\varepsilon$-constraint \\

Makboul et al.~(2022) \cite{Makboul3} 
& \cmark & \cmark & \xmark & \cmark
& OR block assignment
& \xmark
& \xmark
& Block
& ILP + GA + downstream-resource heuristic \\

Aringhieri et al.~(2022) \cite{Aringhieri}
& \cmark & \cmark & \cmark & \cmark
& \xmark
& \xmark
& \xmark
& Block
& ILP + multi-neighborhood local search \\

Deklerck et al.~(2022) \cite{Deklerck} 
& \cmark & \xmark & \xmark & \cmark
& Minimize OR block changes
& \xmark
& \xmark
& Block
& MIP + $\varepsilon$-constraint \\

Kianfar et al.~(2023) \cite{Kianfar} 
& \cmark & \xmark & \xmark & \cmark
& OR block assignment
& \xmark
& \xmark
& Modified block
& MIP + Simulated Annealing + LP \\

Calegari et al.~(2025) \cite{Calegari} 
& \cmark & \xmark & \xmark & \cmark
& \xmark
& \xmark
& \xmark
& Block
& Genetic algorithm + alignment heuristic \\

\midrule

\textbf{This work}
& \cmark & \cmark & \cmark & \cmark
& OR block assignment
& \cmark
& \cmark
& Block/Open/Modified
& MILP + column-generation heuristic \\

\bottomrule
\end{tabular}}
\begin{tablenotes}
\footnotesize
\centering
\item ILP = Integer Linear Program; GA = Genetic Algorithm; LP = Linear Program.
\end{tablenotes}
\end{sidewaystable}

\subsection{Column Generation in OR Scheduling}

The CG approach has been successfully applied to several OR scheduling problems, particularly when compact MILP formulations become computationally challenging. By exploiting problem decomposition, column generation is capable of solving large-scale instances to optimality or near-optimality.

Bargetto et al.~\cite{Bargetto} and Doulabi et al.~\cite{Doulabi} develop branch-and-price-and-cut algorithms for integrated OR scheduling problems involving multiple resources and operational constraints, significantly outperforming compact MILP formulations. Range et al.~\cite{Range} propose a CG framework for patient assignment that explicitly accounts for future arrivals by incorporating stochastic demand information into the pricing problem.

\subsection{Research Gap and Contributions}
\label{section:literature_review:research_gap}

The literature review highlights that only a limited number of studies jointly address the MSSP, the SCAP, and downstream bed management. Moreover, despite the operational importance of reduced weekend staffing, weekend bed occupancy has received little attention in integrated OR scheduling models. This paper addresses these gaps through the following contributions:

\begin{enumerate}
    \item We develop an integrated optimization model that jointly determines the Master Surgical Schedule, assigns surgeons and patients to OR blocks, sequences surgeries, and accounts for bed occupancy during both weekdays and weekends. The model supports block, open, and modified block scheduling policies.

    \item We propose a CG algorithm in which the pricing problem is formulated as a Resource-Constrained Shortest Path Problem (RCSPP), enabling the efficient solution of large-scale instances.

    \item We conduct a comprehensive computational analysis to evaluate: (i) the impact of increasing weekend bed capacity on the number of scheduled patients and overall system performance; (ii) the effect of imposing minimum OR utilization requirements on patient selection; (iii) the relative performance of block, open, and modified block scheduling policies, and (iv) the impact of uncertainty in surgery duration and LoS.
\end{enumerate}

\section{Problem Formulation}
\label{section:problem_formulation}

This section introduces the problem under study. Section~\ref{section:problem_formulation:definition} formally defines the problem, Section~\ref{section:problem_formulation:assumptions} presents the modeling assumptions and derives a structural lemma, and Section~\ref{subsec:math} introduces the MILP formulation.

\subsection{Problem Definition}
\label{section:problem_formulation:definition}

We consider the problem of scheduling elective surgeries in an Operating Theater (OT) consisting of multiple ORs over a multi-day planning horizon while explicitly accounting for surgeon and bed availability. Let $\mathcal{I}$ denote the set of surgeries and $\mathcal{D}$ the set of planning days. Each surgery may either be scheduled or left unscheduled. Each day is divided into shifts, which may consist of a morning and an afternoon shift or a single full-day shift.

ORs are heterogeneous and may differ in terms of specialty compatibility and opening times because of equipment availability or temporary reservations. Let $\mathcal{R}$ denote the set of ORs and $\mathcal{E}$ the set of shifts. Pairing an OR with a shift defines an OR block $(r,e)\in\mathcal{R}\times\mathcal{E}$. Each OR has a fixed daily opening time, which is identical across all shifts.

Each surgery belongs to a specialty, represented by the set $\mathcal{J}$, and is assigned to a surgeon from the corresponding specialty, represented by the set $\mathcal{S}$. Furthermore, each surgery is characterized by a priority score $\pi_i$, an expected duration $\delta_i$, and an expected LoS $\nu_i$.

Patient priority is commonly determined by combining clinical urgency with waiting time, allowing patients with lower initial urgency to progressively gain priority and thereby reducing inequalities in access to surgery \cite{Dery,Mullen,Valente}. Various scoring systems have been proposed in the literature, incorporating procedural, disease-related, and patient-related characteristics \cite{Powers,Rana}. In this work, we assume that priority scores are available and computed according to the hospital's prioritization policy.

Since our objective is to maximize the total priority of scheduled patients, shorter procedures would naturally be favored because they are easier to accommodate within OR blocks and allow a larger number of surgeries to be performed. Consequently, longer procedures would tend to experience longer waiting times. To mitigate this bias, the adopted priority function increases with the expected surgery duration, balancing efficiency and fairness by preventing longer procedures from being systematically postponed. 
The proposed optimization model is independent of the specific priority scoring system adopted. Therefore, the priority score can be replaced by any alternative patient-scoring mechanism or other performance measure, such as expected financial revenue, without requiring structural modifications to the model.

We consider the three scheduling approaches introduced in Section~\ref{section:introduction}:
\begin{enumerate}
\item \textit{Block scheduling}: each OR block is exclusively assigned to a single surgical specialty.
\item \textit{Open scheduling}: OR blocks are shared across specialties and may be used by any compatible surgeon.
\item \textit{Modified block scheduling}: a subset of OR blocks is dedicated to individual specialties, while the remaining blocks are shared across specialties.
\end{enumerate}
In the remainder of the paper, we refer to OR blocks assigned to a single specialty as \emph{dedicated blocks} and to OR blocks accessible by multiple specialties as \emph{shared blocks}. Under block scheduling, all blocks are dedicated; under open scheduling, all blocks are shared; and under modified block scheduling, both dedicated and shared blocks coexist.

Given the OR block configuration, the set of surgical procedures together with their assigned surgeons and specialties, and the available bed capacity, our objective is to determine both the MSS and the detailed operational plan. This requires four interrelated decisions: assigning surgical specialties to OR blocks (\textbf{D1}), assigning surgeons to OR blocks (\textbf{D2}), selecting the surgeries to be scheduled and assigning them to OR blocks (\textbf{D3}), and sequencing surgeries within each OR block (\textbf{D4}).

Decision \textbf{D1} depends on the adopted scheduling policy. In dedicated blocks, exactly one specialty is assigned to each block, whereas shared blocks may accommodate surgeries from multiple specialties. 
Each surgeon can be assigned to at most one OR block per shift, although multiple surgeons may operate within the same block. Each patient belongs to a single specialty and occupies a bed in the corresponding ward throughout the expected LoS. In addition, turnover times between consecutive surgeries are explicitly considered. These may arise from equipment changes, additional cleaning procedures (e.g., following infectious cases), or changes in OR personnel, and depend on both the preceding and the succeeding surgery. 
The objective is to maximize the total priority of scheduled surgeries while minimizing excess weekend bed occupancy and maximizing OR utilization.

The computational complexity of the problem arises from the need to jointly coordinate OR scheduling under limited surgeon and bed availability. The resulting problem generalizes the flow shop scheduling problem and is therefore $\mathcal{NP}$-hard \cite{Ruiz}.

\subsection{Assumptions}
\label{section:problem_formulation:assumptions}

The following assumptions hold throughout the paper:
\textbf{H1)} Surgery durations and patient LoS are deterministic.
\textbf{H2)} Surgeons are pre-assigned to patients.
\textbf{H3)} Scheduled overtime is not allowed.
\textbf{H4)} The waiting list exceeds the planning capacity.
\textbf{H5)} Surgeons cannot move between OR blocks during the same shift, although multiple surgeons may share the same block.
\textbf{H6)} The tactical plan is predefined, including the number of ORs and their opening durations. OR blocks may be half-day or full-day.
\textbf{H7)} Turnover times are deterministic. In particular, \textbf{H7.1}) changing surgeon requires at least as much turnover time as two consecutive surgeries by the same surgeon, and \textbf{H7.2}) turnover times are symmetric.
\textbf{H8)} Emergency patients are not considered.

Assumptions H1--H3 are consistent with those adopted by Bargetto et al.~\cite{Bargetto}. Although surgery durations and LoS are inherently uncertain, they can be estimated using predictive models \cite{Daldossi}. Furthermore, knowledge of the assigned surgeon can improve these estimates \cite{Ibrahim}. Nevertheless, to evaluate the robustness of the proposed approach, Section~\ref{section:analysis:uncertainty} analyzes its performance under uncertain surgery durations.

Assumptions H4--H6 reflect common practice in elective surgery planning. In particular, H5 avoids disruptions caused by surgeons moving between ORs during the same shift, while H6 assumes that strategic decisions regarding OR availability have already been determined. The use of both half-day and full-day OR blocks follows Kianfar et al.~\cite{Kianfar} and increases scheduling flexibility.

Assumption H7 models turnover times as depending only on the pair of consecutive surgeries. Since changing surgeons typically requires additional staff and equipment reconfiguration, turnover times are assumed to be longer whenever consecutive surgeries are performed by different surgeons.

Finally, Assumption H8 is commonly adopted in the literature \cite{Rahimi,Ibrahim,Aringhieri}, reflecting the availability of dedicated resources for emergency patients.

Assumption H7 implies that, once a surgeon starts operating within an OR block, it is never advantageous to interrupt their sequence of surgeries and resume it later. This property is formalized in the following lemma. The proof is provided in~\ref{Appendix2}.

\begin{lemma}
\label{lemma}
Consider a set of patients belonging to the same specialty and potentially assigned to different surgeons. Under Assumptions \textbf{H7.1} and \textbf{H7.2}, any schedule that alternates between surgeons is suboptimal (or, at best, equivalent) to a schedule in which each surgeon performs all of their assigned surgeries consecutively.
\end{lemma}

\subsection{Mathematical Formulation}
\label{subsec:math}
We now present the MILP formulation. Table~\ref{tab:sets_parameters} summarizes the sets, parameters, and decision variables. The constraints are organized into assignment, sequencing, and bed occupancy constraints, followed by the objective function. For notational simplicity, all constraints are written for $r\in\mathcal R$ and $e\in\mathcal E$, but are enforced only for $(r,e)\in\mathcal B$.

\begin{table}[tbp]
\caption{Sets, Parameters, and Decision Variables of the mathematical formulation}
\label{tab:sets_parameters}
\renewcommand{\arraystretch}{1}
\begin{tabular}{lp{14cm}}
\hline
\multicolumn{2}{l}{\textbf{Sets}}\\
\hline
$\mathcal{I}$ & set of patients \\
$\mathcal{J}$ & set of surgical specialties \\
$\mathcal{S}$ & set of surgeons \\
$\mathcal{R}$ & set of ORs \\
$\mathcal{E}$ & set of shifts \\
$\mathcal{D}$ & set of days in the planning horizon \\
$\mathcal{B}\subseteq \mathcal{R}\times \mathcal{E}$ & set of available OR blocks \\
$\mathcal{D}^w \subseteq \mathcal{D}$ & set of weekend days \\
$\mathcal{I}^j \subseteq \mathcal{I}$ & set of patients belonging to specialty $j$ \\
$\mathcal{I}^s \subseteq \mathcal{I}$ & set of patients assigned to surgeon $s$ \\
$\mathcal{J}^r \subseteq \mathcal{J}$ & set of specialties compatible with OR $r$ \\
$\mathcal{E}^d \subseteq \mathcal{E}$ & set of shifts available on day $d$ \\
$\mathcal{B}^b \subseteq \mathcal{B}$ & set of available OR blocks to be assigned to specialties under the modified block scheduling policy \\
\hline
\multicolumn{2}{l}{\textbf{Parameters}}\\
\hline
$\pi_i$ & priority score of patient $i$ \\ 
$\delta_i$ & expected surgery duration of patient $i$ \\ 
$\nu_i$ & expected LoS of patient $i$ \\ 
$\tau_{ii'}$ & turnover time between surgeries of patients $i$ and $i'$ \\ 
$\tilde{\tau}$ & minimum turnover time between two surgeries\\ 
$\sigma_{is}$ & 1 if patient $i$ must be operated on by surgeon $s$, 0 otherwise \\
$\xi_{re}$ & duration (available time) of OR block $(r,e)$ \\
$\beta_j$ & total number of beds available for specialty $j$ \\ 
$\tilde{\beta}_j$ & target bed occupancy for specialty $j$ during weekend days \\
$\hat{\beta}_{jd}$ & number of beds already occupied by previous patients for specialty $j$ on day $d$ \\
$M$ & sufficiently large constant, set equal to twice the maximum shift duration \\ 
\hline
\multicolumn{2}{l}{\textbf{Decision Variables}}\\
\hline
$X_{ire}$ & 1 if patient $i$ is scheduled in OR block $(r,e)$, 0 otherwise \\
$Z_{jre}$ & 1 if specialty $j$ is assigned to OR block $(r,e)$, 0 otherwise \\
$H_{sre}$ & 1 if surgeon $s$ is scheduled in OR block $(r,e)$, 0 otherwise \\
$F_{ire}$ & 1 if patient $i$ is the first patient in OR block $(r,e)$, 0 otherwise \\
$L_{ire}$ & 1 if patient $i$ is the last patient in OR block $(r,e)$, 0 otherwise \\
$O_{ii're}$ & 1 if surgery of patient $i'$ immediately follows that of patient $i$ in block $(r,e)$, 0 otherwise \\
$T_{ire} \ge 0$ & start time of patient $i$’s surgery in block $(r,e)$; 0 if patient $i$ is not scheduled in that block \\
$Y_{id}$ & 1 if patient $i$’s surgery is scheduled on day $d$, 0 otherwise \\
$B_{id}$ & 1 if patient $i$ occupies a bed on day $d$, 0 otherwise \\
$W_{jd} \ge 0$ & number of beds exceeding the target occupancy for specialty $j$ on weekend day $d\in\mathcal{D}^w$\\
\hline
\end{tabular}
\end{table}


\subsubsection{Assignment constraints}
\label{subsec:assignment}
The following constraints define the assignment of specialties, surgeons, and patients to OR blocks.
\begin{align}
\sum_{r\in\mathcal{R}}\sum_{e\in\mathcal{E}} X_{ire}
    &\le 1
&& \forall i\in\mathcal{I}
\label{cns:pat_scheduled_once}\\
T_{ire}
    &\le MX_{ire}
&& \forall i\in\mathcal{I},\,r\in\mathcal{R},\,e\in\mathcal{E}
\label{cns:starttime}\\
X_{ire}
    &\le Z_{j_i re}
&& \forall i\in\mathcal{I},\,r\in\mathcal{R},\,e\in\mathcal{E}
\label{cns:spec_pat}\\
X_{ire}
    &\le H_{s_i re}
&& \forall i\in\mathcal{I},\,r\in\mathcal{R},\,e\in\mathcal{E}
\label{cns:sur_pat}\\
\sum_{j\in\mathcal{J}^r} Z_{jre}
    &\le 1
&& \forall r\in\mathcal{R},\,e\in\mathcal{E}
\label{cns:max1speca}\\
\sum_{j\in\mathcal{J}\setminus\mathcal{J}^r} Z_{jre}
    &=0
&& \forall r\in\mathcal{R},\,e\in\mathcal{E}
\label{cns:max1specb}\\
\sum_{r\in\mathcal{R}} H_{sre}
    &\le 1
&& \forall s\in\mathcal{S},\,e\in\mathcal{E}
\label{cns:spec_sura}\\
H_{sre}
    &\le \sum_{i\in\mathcal{I}^s}X_{ire}
&& \forall s\in\mathcal{S},\,r\in\mathcal{R},\,e\in\mathcal{E}
\label{cns:spec_surb}\\
T_{ire}+\delta_i
    &\le \xi_{re}+M(1-X_{ire})
&& \forall i\in\mathcal{I},\,r\in\mathcal{R},\,e\in\mathcal{E}.
\label{cns:blocklength}
\end{align}

Constraints~\eqref{cns:pat_scheduled_once} ensure that each patient is scheduled in at most one OR block over the planning horizon, while constraints~\eqref{cns:starttime} set the surgery start time to zero whenever patient $i$ is not assigned to OR block $(r,e)$.
Constraints~\eqref{cns:spec_pat}--\eqref{cns:sur_pat} require that a patient can be assigned to an OR block only if both the corresponding specialty and surgeon are assigned to that block.
Constraints~\eqref{cns:max1speca}--\eqref{cns:max1specb} ensure that each dedicated OR block is assigned to at most one compatible specialty. These constraints are enforced only for dedicated blocks and are omitted for shared blocks.
Constraints~\eqref{cns:spec_sura}--\eqref{cns:spec_surb} impose that each surgeon can be assigned to at most one OR block per shift and only if at least one of their patients is scheduled in that block. 
Finally, Constraint~\eqref{cns:blocklength} ensures that every scheduled surgery is completed within the available duration of its assigned OR block.

\subsubsection{Ordering Constraints}
\label{subsec:ordering}

The following constraints determine the sequencing of surgeries within each OR block.
\begin{align}
T_{i're} &\ge T_{ire}+\delta_i+\tau_{ii'}-M(1-O_{ii're})
&& \forall i\neq i'\in\mathcal I,\; r\in\mathcal R,\; e\in\mathcal E
\label{cns:ordering_start}\\
\sum_{i\in\mathcal I}F_{ire}
&\le 1
&& \forall r\in\mathcal R,\;e\in\mathcal E
\label{cns:firsta}\\
\sum_{i\in\mathcal I}L_{ire}
&\le 1
&& \forall r\in\mathcal R,\;e\in\mathcal E
\label{cns:lasta}\\
F_{ire}
&\le X_{ire}
&& \forall i\in\mathcal I,\;r\in\mathcal R,\;e\in\mathcal E
\label{cns:firstb}\\
L_{ire}
&\le X_{ire}
&& \forall i\in\mathcal I,\;r\in\mathcal R,\;e\in\mathcal E
\label{cns:lastb}\\
O_{i'ire}
&\le X_{ire}
&& \forall i\in\mathcal I,\;r\in\mathcal R,\;e\in\mathcal E
\label{cns:pred}\\
O_{ii're}
&\le X_{ire}
&& \forall i\in\mathcal I,\;r\in\mathcal R,\;e\in\mathcal E
\label{cns:succ}\\
X_{ire}
&=
F_{ire}+\sum_{i'\in\mathcal I,\;i'\neq i}O_{i'ire}
&& \forall i\in\mathcal I,\;r\in\mathcal R,\;e\in\mathcal E
\label{cns:first}\\
X_{ire}
&=
L_{ire}+\sum_{i'\in\mathcal I,\;i'\neq i}O_{ii're}
&& \forall i\in\mathcal I,\;r\in\mathcal R,\;e\in\mathcal E
\label{cns:last}\\
T_{ire}
&\le (1-F_{ire})M
&& \forall i\in\mathcal I,\;r\in\mathcal R,\;e\in\mathcal E
\label{cns:firsttime}
\end{align}

Constraint~\eqref{cns:ordering_start} defines surgery start times and prevents overlapping surgeries within the same OR block.
Constraints~\eqref{cns:firsta}--\eqref{cns:lasta} ensure that each OR block contains at most one first patient and one last patient.
Constraints~\eqref{cns:firstb}--\eqref{cns:succ} enforce consistency between the sequencing variables and the assignment variables, ensuring that only scheduled patients can appear as first, last, predecessor, or successor in an OR block.
Constraints~\eqref{cns:first} and~\eqref{cns:last} ensure that every scheduled patient has exactly one predecessor or is the first patient in the sequence, and exactly one successor or is the last patient in the sequence.
Finally, constraints~\eqref{cns:firsttime} sets the start time of the first surgery in each OR block to zero.

\subsubsection{Bed occupancy constraints}
\label{subsec:bed}

The following constraints define patient bed occupancy and enforce capacity limits.
\begin{align}
Y_{id}
    &= \sum_{e\in\mathcal{E}^d}\sum_{r\in\mathcal{R}} X_{ire}
&& \forall i\in\mathcal{I},\, d\in\mathcal{D}
\label{cns:surgeryday}\\
\sum_{d'=d}^{\min\{d+\nu_i,|\mathcal{D}|\}} B_{id'}
    &\ge
    \min\{\nu_i,|\mathcal{D}|-d\}Y_{id}
&& \forall i\in\mathcal{I},\, d\in\mathcal{D}
\label{cns:bedoccupancy}\\
\sum_{i\in\mathcal{I}^j} B_{id}
    &\le
    \beta_j-\hat{\beta}_{jd}
&& \forall j\in\mathcal{J},\, d\in\mathcal{D}
\label{cns:bedcapacity}\\
W_{jd}
    &\ge
    \sum_{i\in\mathcal{I}^j}B_{id}-\tilde{\beta}_j+\hat{\beta}_{jd}
&& \forall j\in\mathcal{J},\, d\in\mathcal{D}^w
\label{cns:weekendbedsa}\\
W_{jd}
    &\ge 0
&& \forall j\in\mathcal{J},\, d\in\mathcal{D}^w.
\label{cns:weekendbedsb}
\end{align}

Constraint~\eqref{cns:surgeryday} defines the binary variable $Y_{id}$, indicating whether patient $i$ undergoes surgery on day $d$.
Constraint~\eqref{cns:bedoccupancy} ensures that, once operated on, a patient occupies a bed for the entire expected LoS.
Constraint~\eqref{cns:bedcapacity} enforces the bed-capacity limits for each specialty.
Finally, Constraints~\eqref{cns:weekendbedsa}--\eqref{cns:weekendbedsb} define the excess weekend bed occupancy variable $W_{jd}$, which measures the number of beds exceeding the target occupancy level for specialty $j$ on weekend day $d$.

\subsubsection{Objective Function}
\label{sub:obj}
The model considers three objectives: (i) maximizing the total priority of scheduled patients, (ii) minimizing excess weekend bed occupancy, and (iii) maximizing OR utilization. OR utilization is computed only over blocks with at least one scheduled surgery and includes the minimum turnover time between consecutive surgeries.
\begin{align}
\label{eq:obj}
    \max  \quad & \sum_{e\in\mathcal{E}}\sum_{r\in\mathcal{R}}\sum_{i\in\mathcal{I}} \pi_{i} X_{ire} \\
    \min \quad & \sum_{j\in\mathcal{J}}\sum_{d\in\mathcal{D}^w} W_{jd} \\
    \max \quad &\sum_{r\in\mathcal{R}}\sum_{e\in\mathcal{E}|\sum_{i\in\mathcal{I}}X_{ire}>0}\frac{\sum_{i\in\mathcal{I}}\left[X_{ire}(\delta_i+\tilde{\tau})\right] -\tilde{\tau}}{\xi_{re}}
\end{align}

Rather than optimizing the three objectives simultaneously, patient priority is treated as the primary objective, while excess weekend bed occupancy and OR utilization are controlled through constraints. Let $\Theta_{1_{jd}}$ denote the maximum allowable excess weekend occupancy for specialty $j$ on day $d$, and let $\Theta_2$ denote the minimum required OR utilization. The resulting optimization model is
\begin{align}
    \label{obj:def}
    \max \quad & \sum_{e\in\mathcal{E}}\sum_{r\in\mathcal{R}}\sum_{i\in\mathcal{I}} \pi_{i} X_{ire} & \\
    \text{s.t.} \quad  &  \eqref{cns:pat_scheduled_once}-\eqref{cns:weekendbedsb} \\
    &  W_{jd}\le \Theta_{1_{jd}} & \quad \forall j\in\mathcal{J}, d\in\mathcal{D}^w \\
    & \sum_{r\in\mathcal{R}}\sum_{e\in\mathcal{E} | \sum_{i\in\mathcal{I}}X_{ire} > 0}\frac{\sum_{i\in\mathcal{I}}\left[X_{ire}(\delta_i+\tilde{\tau})\right]-\tilde{\tau}}{\xi_{re}} \ge \Theta_2 &
\end{align}
\section{Column Generation Formulation}
\label{section:cg}
In this section, we present the column generation approach used to solve the problem. In Section \ref{section:cg:overview}, we provide a brief overview of the CG procedure used in our framework. In Section \ref{section:cg:mp}, we present the mathematical formulation of the master problem. In Section \ref{section:cg:sp}, we derive the dual problem, formulate the pricing problem, and present the RCSPP formulation of the subproblem, while the MILP formulation is provided in the Supplementary Material. Finally, in Section \ref{section:cg:algorithm}, we describe the implementation of the column generation algorithm along with its pseudocode.

\subsection{Column Generation Overview}
\label{section:cg:overview}
Column generation is a solution approach based on the primal simplex algorithm. Instead of explicitly considering all variables, it iteratively searches for variables with negative reduced cost, or proves that no such variables exist. For further details on the column generation methodology, we refer the reader to the standard literature, such as \cite{Desrosiers, Uchoa}. 

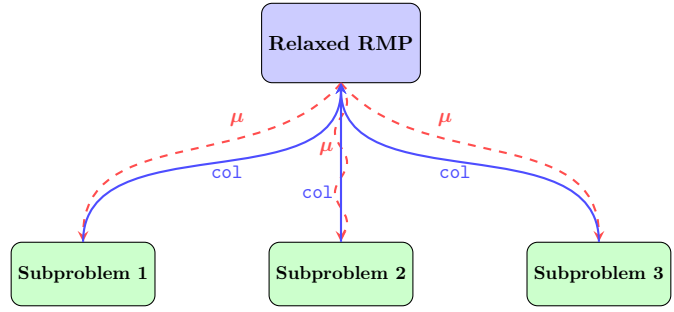
\begin{figure}[tb]
\caption{Framework of the CG procedure.}
    \label{figure:cg}
\tikzset{
    master/.style={rectangle, rounded corners, minimum width=3cm, minimum height=1.5cm, 
                   text centered, draw=black, fill=blue!20, font=\bfseries},
    subproblem/.style={rectangle, rounded corners, minimum width=2.5cm, minimum height=1.2cm,
                       text centered, draw=black, fill=green!20, font=\small\bfseries},
    arrow/.style={thick,->,>=stealth},
    dual/.style={arrow, color=red!70, dashed},
    column/.style={arrow, color=blue!70}
}
\begin{center}
\begin{tikzpicture}[scale=0.7, transform shape, node distance=2cm]
    \node (master) [master] {Relaxed RMP};
    
    \node (sub1) [subproblem, below left=3cm and 2cm of master] {Subproblem 1};
    \node (sub2) [subproblem, below=3cm of master] {Subproblem 2};
    \node (sub3) [subproblem, below right=3cm and 2cm of master] {Subproblem 3};
    
    \draw [dual] (master.south) to[out=-130, in=90] node[left, pos=0.3, yshift=3mm, text=red!70] {$\boldsymbol{\mu}$} (sub1.north);
    \draw[dual,
      decorate,
      decoration={snake, amplitude=0.8mm, segment length=8mm}]
    (master.south)
    to[out=-90,in=90] node[left, pos=0.3, yshift=-3mm, text=red!70] {$\boldsymbol{\mu}$} 
    (sub2.north);
    \draw [dual] (master.south) to[out=-50, in=90] node[right, pos=0.3, yshift=3mm, text=red!70] {$\boldsymbol{\mu}$} (sub3.north);
    
    \draw [column] (sub1.north) to[out=90, in=-90] node[left, pos=0.6, yshift=-3mm, text=blue!70] {\texttt{col}} (master.south);
    \draw [column] (sub2.north) to[out=90, in=-90] node[left, pos=0.3, text=blue!70] {\texttt{col}} (master.south);
    \draw [column] (sub3.north) to[out=90, in=-90] node[right, pos=0.6, yshift=-3mm, text=blue!70] {\texttt{col}} (master.south);
\end{tikzpicture}
\end{center}
\end{figure}

The CG procedure used in our framework is illustrated in Figure \ref{figure:cg} and in Algorithm \ref{alg:cg_steps}.
We iteratively solve the relaxed Reduced Master Problem (RMP), compute the dual values $\mu$, and use them to define and solve the subproblems, thereby generating new columns \texttt{col}. 
When no column with positive reduced cost can be found, then all columns required for the optimal solution of the relaxed RMP have been generated, and the problem can be solved to optimality. Otherwise, a 
solution with positive reduced cost is identified, and the corresponding column is added to the relaxed RMP. The relaxed RMP is then re-optimized, and this process is repeated until no column with positive reduced cost can be found.

\begin{algorithm}[tb]
\caption{Column Generation procedure}
\label{alg:cg_steps}
\begin{algorithmic}[1]
\Statex \textbf{Input:} $\mathcal{R}$ (OR rooms), $\mathcal{E}$ (shifts), $\xi_{re}$ (opening times), $\mathcal{I}$ (patients), $\mathcal{S}$ (surgeons), $\mathcal{J}$ (specialties)
\Statex \textbf{Output:} Integer solution of the RMP
\State Generate an initial feasible set of columns.
\State Solve the RMP relaxation and obtain dual variables $\boldsymbol{\mu}$.
\State Solve the subproblems and generate columns \texttt{col} with positive reduced cost.
\State \textbf{while} at least one column with positive reduced cost exists \textbf{do}
\State \hspace{1em} Add generated columns to the relaxed RMP.
\State \hspace{1em} Re-optimize the relaxed RMP and update dual variables $\boldsymbol{\mu}$.
\State \hspace{1em} Resolve the subproblems.
\State \textbf{end while}
\State Solve the integer RMP using all generated columns.
\State \Return Optimal integer RMP solution
\end{algorithmic}
\end{algorithm}

We decompose the problem as follows.
\begin{description}
\item[\textbf{Master Problem}:] It takes as input a set of feasible schedules and determines the optimal combination of schedules that maximizes the objective function~\eqref{obj:def}, while enforcing constraints that govern interactions between different OR blocks. These include surgeon availability (a surgeon cannot be assigned to more than one OR in the same shift), bed capacity constraints, and minimum OR utilization requirements. Although the original decision variables are binary, at each iteration of the column generation algorithm we solve a Relaxed RMP where the variables are continuous and take values in the interval $[0,1]$.
\item[\textbf{Subproblem}:] It restricts the search for improving columns to a single OR block and identifies the schedule that maximizes the reduced cost for that block. Since a MILP formulation of the subproblem is computationally expensive, we propose a RCSPP formulation. 
\end{description}




\subsection{Master Problem}
\label{section:cg:mp}
Let $\Omega_{re}$ be the set of all feasible surgery schedules for the OR block $(r,e)$, and let $k \in \Omega_{re}$ denote a feasible schedule.
A schedule is represented by a set of ordered scheduled surgeries, together with the set of surgeons assigned to the OR block, which corresponds to the surgeons required for the surgeries. A schedule is feasible if:
(i) the sum of surgery durations and turnover times does not exceed $\xi_{re}$;
(ii) surgeries do not overlap;
(iii) no patient is scheduled more than once;
(iv) all patients in the schedule belong to the same specialty (only for the dedicated OR blocks).

Given a schedule, let $\alpha_{ik}$ be equal to 1 if patient $i$ is included in schedule $k$ and 0 otherwise; let $\epsilon_k=\sum_{i\in\mathcal I}\pi_i\alpha_{ik}$ denote the contribution of schedule $k$ to the objective function; let $\rho_{sk}$ be equal to 1 if surgeon $s$ is selected in schedule $k$ and 0 otherwise; and let $\psi_{jdk}\ge 0$ denote the number of beds occupied by schedule $k$ for specialty $j$ on day $d$.

We must determine which schedules to select in order to maximize the objective function while satisfying all the imposed constraints. The decision variable is denoted by $C_k$ and takes the value $1$ if schedule $k$ is chosen, and $0$ otherwise. Thus, we aim to solve
\begin{align}
    \label{obj:master}
    \max &\quad \sum_{r\in \mathcal{R}}\sum_{e\in \mathcal{E}}\sum_{k\in\Omega_{re}} \epsilon_{k} C_{k} \\
    \text{s.t. } \sum_{r\in \mathcal{R}}\sum_{e\in \mathcal{E}}\sum_{k\in\Omega_{re}} \alpha_{ik}C_{k} &\le 1 & \forall i\in \mathcal{I} & \quad [\mu^1_i \ge 0] \label{cns:patientonce}\\
    \sum_{r\in \mathcal{R}} \sum_{k\in \Omega_{re}} \rho_{sk}C_k &\le 1 & \forall s\in \mathcal{S}, e\in \mathcal{E} & \quad [\mu^2_{se} \ge 0] \label{cns:surgeonone} \\
    \sum_{k\in\Omega_{re}} C_{k} &\le 1 & \forall r\in \mathcal{R}, e\in \mathcal{E} & \quad [\mu^3_{re} \ge 0] \label{cns:scheduleone} \\
    \sum_{r\in \mathcal{R}}\sum_{e\in \mathcal{E}} \sum_{k\in \Omega_{re}} \psi_{jdk}C_k &\le \beta_{j} - \hat{\beta}_{jd} & \forall j\in \mathcal{J}, d\in \mathcal{D} & \quad [\mu^4_{jd} \ge 0] \label{cns:beds} \\
    \sum_{r\in \mathcal{R}}\sum_{e\in \mathcal{E}}\sum_{k\in \Omega_{re}} \psi_{jdk}C_k &\le \tilde{\beta}_{j} - \hat{\beta}_{jd}+ \Theta_{1_{jd}} & \forall j\in \mathcal{J}, d\in \mathcal{D}^w & \quad [\mu^5_{jd} \ge 0] \label{cns:bedsweekend}
    \\
    \sum_{r\in \mathcal{R}}\sum_{e\in \mathcal{E}}\sum_{k\in\Omega_{re}}C_k \left\{\sum_{i\in \mathcal{I}} \left[\alpha_{ik}(\delta_i+\tilde{\tau})\right]  -\tilde{\tau}\right\}&\ge \Theta_2\sum_{r\in \mathcal{R}}\sum_{e\in \mathcal{E}}\xi_{re} & & \quad [\mu^6 \le 0] \label{cns:utilization}
\end{align}

Constraints \eqref{cns:patientonce} ensure that each patient is scheduled at most once. Constraints \eqref{cns:surgeonone} ensure that each surgeon works in at most one OR block at a time. Constraints \eqref{cns:scheduleone} ensure that at most one schedule is selected for each OR block. Constraints \eqref{cns:beds} enforce the bed capacity limit during weekdays, while constraints \eqref{cns:bedsweekend} enforce the bed capacity limit during weekends, allowing for $\Theta_{1_{jd}}$ additional beds.
Finally, constraints \eqref{cns:utilization} ensure that the total utilization is at least $\Theta_2$ of the total capacity of the OR blocks.

In each column generation iteration, we solve the relaxed RMP, i.e., model \eqref{obj:master}–\eqref{cns:utilization} with $C_k \ge 0$ rather than $C_k\in\{0,1\}$; the bound $C_k \le 1$ is naturally enforced by constraints \eqref{cns:scheduleone}.


A large value of $\Theta_2$ may render the initial RMP infeasible under the columns generated at initialization. We therefore modify the RMP by introducing a nonnegative slack variable $\boldsymbol{\mathcal{S}} \geq 0$, which guarantees feasibility of the initial RMP. The slack variable absorbs any initial violation of the constraint and allows the algorithm to proceed until sufficient columns are generated to recover feasibility without reliance on the slack. This follows the usual column-generation strategy of enforcing initial feasibility via auxiliary variables in the master problem.
Constraints~\eqref{cns:utilization} become
\begin{equation}
\sum_{r\in \mathcal{R}}\sum_{e\in \mathcal{E}}\sum_{k\in\Omega_{re}}C_k \left\{\sum_{i\in \mathcal{I}} \left[\alpha_{ik}(\delta_i+\tilde{\tau})\right]  -\tilde{\tau}\right\} + \boldsymbol{\mathcal{S}} \ge \Theta_2\sum_{r\in \mathcal{R}}\sum_{e\in \mathcal{E}}\xi_{re}
\end{equation}

and the objective function~\eqref{obj:master} becomes
\begin{equation}
\label{obj:master2}
\max \quad \sum_{r\in \mathcal{R}}\sum_{e\in \mathcal{E}}\sum_{k\in\Omega_{re}} \epsilon_{k} C_{k} - \boldsymbol{\mathcal{S}}.
\end{equation}

When column generation terminates, the integer RMP is solved without the variable $\boldsymbol{\mathcal{S}}$.

\subsection{Subroblem}
\label{section:cg:sp}
Given the RMP, the corresponding dual problem is defined as follows:
\begin{align}
    \min \quad  & \sum_{i\in \mathcal{I}}\mu^1_i + 
    \sum_{s\in \mathcal{S}}\sum_{e\in \mathcal{E}} \mu^2_{se} + 
    \sum_{r\in \mathcal{R}}\sum_{e\in \mathcal{E}} \mu^3_{re} + 
    \sum_{j\in \mathcal{J}}\sum_{d\in \mathcal{D}} (\beta_{j}-\hat{\beta}_{jd})\mu^4_{jd} + \\
    & \quad \sum_{j\in \mathcal{J}}\sum_{d\in \mathcal{D}^w} (\tilde{\beta}_j-\hat{\beta}_{jd}+\Theta_{1_j}) \mu^5_{jd} + \Theta_2  \sum_{r\in \mathcal{R}}\sum_{e\in \mathcal{E}} \xi_{re}\mu^6 \notag \\
    \text{s.t.} \quad &
\sum_{i\in \mathcal{I}} \alpha_{ik} \mu^1_i
+ \sum_{s\in \mathcal{S}} \rho_{sk} \mu^2_{se}
+ \mu^3_{re}
+ \sum_{j\in \mathcal{J}} \sum_{d\in \mathcal{D}} \psi_{jdk} \mu^4_{jd}
+ \sum_{j\in \mathcal{J}} \sum_{d\in \mathcal{D}^w} \psi_{jdk} \mu^5_{jd}  \\
&+ \left(\sum_{i\in \mathcal{I}} \alpha_{ik}(\delta_i+\tilde{\tau}) - \tilde{\tau}\right)\mu^6
\ge \epsilon_k
\quad \forall r\in \mathcal{R}, e\in \mathcal{E}, k\in \Omega_{re} \notag
\\
&\mu^1_i \ge 0 \ \forall i\in \mathcal{I}, \quad
\mu^2_{se} \ge 0 \ \forall s\in \mathcal{S}, e\in \mathcal{E}, \quad
\mu^3_{re} \ge 0 \ \forall r\in \mathcal{R}, e\in \mathcal{E}, \notag \\
&\mu^4_{jd} \ge 0 \ \forall j\in \mathcal{J}, d\in \mathcal{D}, \quad
\mu^5_{jd} \ge 0 \ \forall j\in \mathcal{J}, d\in \mathcal{D}^w, \quad
\mu^6 \le 0
\end{align}

Thus, for each block $(r,e)$, the pricing problem is: 
\begin{align}
\max\; \hat{\epsilon}_k 
= \max \Bigg\{
&\sum_{i\in \mathcal{I}}\pi_i\alpha_{ik} -
\Big(
\sum_{i\in \mathcal{I}} \alpha_{ik} \mu^1_i 
+ \sum_{s\in \mathcal{S}}\rho_{sk} \mu^2_{se} 
+ \mu^3_{re} \\
&\quad
+ \sum_{j\in \mathcal{J}}\sum_{d'\in\mathcal{D}}\psi_{jd'k} \mu^4_{jd'} 
+ \sum_{j\in \mathcal{J}}\sum_{d'\in\mathcal{D}^w} \psi_{jd'k} \mu^5_{jd'} 
+ \left\{\sum_{i\in \mathcal{I}}\left[\alpha_{ik}(\delta_i+\tilde{\tau})\right]-\tilde{\tau}\right\} \mu^6
\Big)
\Bigg\} \notag
\end{align}

For dedicated OR blocks, the structure of the problem allows the subproblem to be decomposed by specialty. Specifically, instead of solving a single subproblem for each OR block that considers all specialties simultaneously, we fix one specialty at a time and solve a separate subproblem for each OR block–specialty pair. This decomposition increases the number of generated columns while reducing the number of patients considered in each subproblem, thereby improving computational efficiency.


We formulate the subproblem as an RCSPP. Starting from a dummy node \textit{start} and ending at a dummy node \textit{end}, we define a graph where nodes represent patients. Arcs are associated with a profit corresponding to the gain from scheduling a patient and a resource consumption equal to the surgical time, including turnover.
We denote the graph as $G = (V, A)$, where $V = \{\textit{start}, i_1, i_2, \ldots, i_{|\mathcal{I}|}, \textit{end}\}$ is the set of nodes, and $A \subseteq V \times V$ is the set of arcs connecting the nodes.

A path from \textit{start} to \textit{end} is feasible if:
\begin{enumerate}
\item Its total duration does not exceed $\xi_{re}$:
\[
\sum_{(i,i') \in \text{path}} (\delta_i + \tau_{ii'}) \le \xi_{re}.
\]
\item The path is elementary, i.e., each node is visited at most once.
\item For dedicated OR blocks, all scheduled patients belong to the same specialty.
\end{enumerate}

Again, due to the structure of the problem, we can omit the third condition and restrict the set of patients to a single specialty, yielding $V = \{\textit{start}, i_1, i_2, \ldots, i_{|\mathcal{I}^j|}, \textit{end}\}$, and solve each subproblem for each couple OR block-specialty.

We aim to find the path from \textit{start} to \textit{end} that maximizes the following expression:

\begin{equation}
    \max_{\text{path}} \sum_{(i,i') \in \text{path}} \tilde{\pi}_{i'} - \sum_{s\in \mathcal{S}}\rho_s\mu^2_{s} - \mu^3  +\tilde{\tau}\mu^6,
\end{equation}
where $d$ is such that $e\in \mathcal{E}^d$, and $\tilde{\pi}_{i'}$ is defined as follows:
\begin{equation}
    \tilde{\pi}_{i'} = \pi_{i'} -
    \mu^1_{i'} - 
    \sum_{d'=d}^{\min\{d+\nu_i', |\mathcal{D}|\}}\mu^4_{j_i' d'} - 
    \sum_{d'=d, d'\in\mathcal{D}^w}^{\min\{d+\nu_i', |\mathcal{D}|\}}\mu^5_{j_i' d'} - (\delta_{i'}+\tilde{\tau})\mu^6.
\end{equation}

An example of a graph with four patients and two surgeons is shown in Figure \ref{Grafo_shift}. By Lemma~\ref{lemma}, the graph can be constructed in this way without counting the surgeon’s dual variable twice, since switching surgeons within a schedule is never optimal due to the additional time required. 
\begin{figure}[tb]
\caption{Example of a graph with four patients and two surgeons. \textit{Start} and \textit{End} denote the dummy source and sink nodes, respectively. Nodes $i_1$ and $i_2$ represent the patients assigned to surgeon $s_1$, while nodes $i_3$ and $i_4$ represent the patients assigned to surgeon $s_2$. All arcs of the graph are shown, but gains (in orange) and resource consumption (in blue) are reported only for a subset of them. In the implementation, since the RCSPP is typically formulated as a cost minimization problem, we solve it by minimizing the costs multiplied by $-1$, which is equivalent to maximizing the gains.}
\label{Grafo_shift}
    \centering
\begin{tikzpicture}[scale=0.9, transform shape, 
    node distance=2cm and 2cm,
    every node/.style={font=\sffamily},
    patient/.style={circle, draw, minimum size=10mm, inner sep=0pt},
    box/.style={rectangle, draw, minimum width=12mm, minimum height=12mm, align=center},
    arrow/.style={<->, very thin, shorten >=2pt, shorten <=2pt}
]

\node[box] (start) at (0,0) {Start};
\node[box] (end) at (8,0) {End};

\node[patient] (p1) at (2.5, 1.5) {$i_1$};
\node[patient] (p2) at (2.5, -1.5) {$i_2$};

\node[patient] (p3) at (5.5, 1.8) {$i_3$};
\node[patient] (p4) at (5.5, -1.8) {$i_4$};

\foreach \a in {1,2}
  \foreach \b in {3,4}
    \draw[arrow] (p\a) -- (p\b);
\draw[arrow] (p1) -- (p2);
\draw[arrow] (p3) -- (p4);

\foreach \i in {1,2,3,4}
  \draw[->, very thin] (start) -- (p\i);
\foreach \i in {1,2,3,4}
  \draw[->, very thin] (p\i) -- (end);

\coordinate (mid12) at ($(p1)!0.5!(p2)$);
\draw[dashed] (mid12) ellipse [x radius=1.2cm, y radius=2.2cm];
\node[above=2.2cm of mid12] {$s_1$};

\coordinate (mid34) at ($(p3)!0.5!(p4)$);
\draw[dashed] (mid34) ellipse [x radius=1.2cm, y radius=2.4cm];
\node[above=2.4cm of mid34] {$s_2$};

\draw[->, thick, orange] 
    (start) -- 
    node[above, sloped, fill=white, inner sep=1pt, font=\footnotesize]{$\tilde{\pi}_4-\mu^2_{s_2}$} 
    node[below, sloped, fill=white, inner sep=1pt, font=\small]{\textcolor{blue}{$\delta_4$}} 
    (p4);
    
\draw[->, thick, orange] 
    (p4) -- 
    node[above, sloped, fill=white, inner sep=1pt, font=\footnotesize]{$\tilde{\pi}_3$} 
    node[below, sloped, fill=white, inner sep=1pt, font=\small]{\textcolor{blue}{$\delta_3+\tau_{43}$}} 
    (p3);
    
\draw[->, thick, orange] 
    (p1) -- 
    node[above, sloped, fill=white, inner sep=1pt, font=\footnotesize]{$\tilde{\pi}_4-\mu^2_{s_2}$} 
    node[below, sloped, fill=white, inner sep=1pt, font=\small]{\textcolor{blue}{$\delta_4 + \tau_{14}$}} 
    (p4);

\draw[->, thick, orange] 
    (p3) -- 
    node[above, sloped, fill=white, inner sep=1pt, font=\footnotesize]{$-\mu^3+\tilde{\tau}\mu^6$} 
    (end);
\end{tikzpicture}
\end{figure}

\subsection{Column Generation Algorithm}
\label{section:cg:algorithm}
In this section, we describe the steps implemented in the column generation algorithm.

\begin{description}
    \item[\textbf{Input:}] ORs, shifts, opening times of each OR, and the set of patients, surgeons, and specialties, along with all their associated information.

    \item [\textbf{Initial Columns Construction:}]
We first generate an initial set of columns ensuring the feasibility of the master problem. Patients are grouped by surgeon and sorted by priority code. For each shift, surgeons are assigned to available ORs (leaving some ORs empty if surgeons are insufficient). Patients are then sequentially assigned to their surgeon's OR blocks according to priority, with start times computed by adding surgery and turnover times. If a patient cannot fit within the current OR block, the next block is considered. The procedure terminates when all OR blocks or patients are assigned. Finally, single-patient schedules are generated as dummy columns in the CG framework by assigning each patient to every compatible OR block.

    \item [\textbf{Relaxed RMP:}] We solve the Relaxed RMP to obtain the dual variable values.

   \item [\textbf{Subproblem:}]
According to the formulation, we solve a subproblem for each pair $(r,e)$ in shared OR blocks, and for each pair $(r,e)$ and specialty $j$ in dedicated OR blocks, to determine the optimal specialty and patient set given the dual variables.
For shared ORs operating two shifts, the subproblem is solved only for the morning shift and the resulting solution is applied to the evening shift, since the two shifts are equivalent (i.e., same day and opening time). Similarly, for dedicated ORs, we solve one subproblem for each day, opening time, and specialty, and apply the solution to all compatible OR blocks with the same configuration.
Therefore, if $|\mathcal{R}_2|$ and $|\mathcal{R}_1|$ denote the number of ORs operating two and one shifts, respectively, the number of subproblems for shared OR blocks is reduced from $|\mathcal{R}_2|\times|\mathcal{E}|+|\mathcal{R}_1|\times|\mathcal{D}|$ to $(|\mathcal{R}_2|+|\mathcal{R}_1|)\times|\mathcal{D}|$. For dedicated OR blocks, the number of subproblems is reduced from $|\mathcal{R}_2|\times|\mathcal{E}|\times|\mathcal{J}|+|\mathcal{R}_1|\times|\mathcal{D}|\times|\mathcal{J}|$ to $R_{\text{types}}\times|\mathcal{D}|\times|\mathcal{J}|$, where $R_{\text{types}}$ denotes the number of distinct OR opening configurations (equal to 2 in our setting).
New columns are generated by solving the RCSPP formulation.
\item [\textbf{Final Solution:}]
After the CG process, we retain all columns generated during the solution of the relaxed RMP. The enriched column set is then used to solve the RMP with integrality constraints, requiring $C_k$ to be binary. This provides a feasible solution, although optimality is not guaranteed without a full branch-and-price procedure. The proposed approach achieves high-quality solutions while substantially reducing computational effort compared with a complete branch-and-price scheme.
\end{description}

It is worth noting that, for the RCSPP formulation, the problem is initially solved heuristically. If a column with positive reduced cost is identified, the algorithm adds the column to the RMP and proceeds to the next iteration. Otherwise, the problem is solved exactly. If no column with positive reduced cost is found, the algorithm moves on to the next subproblem; otherwise, the column is added to the RMP and the heuristic procedure is restarted.

The pseudocode for obtaining an initial feasible schedule is presented in Algorithm~\ref{alg:cg}, while the CG loop is described in Algorithm~\ref{alg:cgalgo}.

\begin{algorithm}[tb]
\caption{Phase 1: Build initial feasible solution}
\label{alg:cg}
\begin{algorithmic}[1]

\Statex \textbf{Input:} $\mathcal{R}$ (OR rooms), $\mathcal{E}$ (shifts), $\xi_{re}$ (opening times),
                        $\mathcal{I}$ (patients), $\mathcal{S}$ (surgeons), $\mathcal{J}$ (specialties)
\Statex \textbf{Output:} $\textsc{Initial Schedule}$

\For{each surgeon $s \in \mathcal{S}$}
    \State $\mathtt{surgeon\_patients}[s] \leftarrow \textsc{SortByPriority}(\mathcal{I}^s)$  \Comment{assign sorted patients to surgeons}
\EndFor
\For{each shift $e\in \mathcal{E}$}
    \For{each OR $r\in\mathcal{R}$}
        \State $s \leftarrow \textsc{AssignSurgeon}(r, e)$
        \For{each patient $i \in \mathtt{surgeon\_patients}[s]$}
            \If{$i$ fits in remaining time of block $(r,e)$}
                \State Schedule $i$ in $\textsc{Initial Schedule}$; set $i.\mathtt{start} \leftarrow \textsc{ComputeStartTime}(i)$
            \Else
                \State Move to next available OR
            \EndIf
        \EndFor
    \EndFor
\EndFor
\For{each patient $i \in \mathcal{I}$}
    \For{each shift $e \in \mathcal{E}$}
        \For{each OR $r\in\mathcal{R}$}
            \If{OR $r$ compatible with the specialty of $i$ and $i$ fits in block $(r, e)$}
                \State Schedule $i$ in $\textsc{Initial Schedule}$; set $i.\mathtt{start} \leftarrow 0$
            \Else
                \State Move to next available OR
             \EndIf
        \EndFor
    \EndFor
 \EndFor
\State \Return $\textsc{Initial Schedule}$

\end{algorithmic}
\end{algorithm}

\begin{algorithm}[tb]
\caption{Phase 2 — Column generation loop}
\label{alg:cgalgo}
\begin{algorithmic}[1]

\Statex \textbf{Input:} $\mathcal{R}$ (OR rooms), $\mathcal{E}$ (shifts), $\xi_{re}$ (opening times),
                        $\mathcal{I}$ (patients), $\mathcal{S}$ (surgeons), $\mathcal{J}$ (specialties), $\textsc{Initial Schedule}$
\Statex \textbf{Output:} Optimal solution of the integer RMP

\Repeat 
\Until{time limit exceeded \textbf{or} maximum iterations reached \textbf{or} $\mathtt{improved} = \texttt{False}$}
    \State $\boldsymbol{\mu} \leftarrow \textsc{SolveRelaxedRMP}()$ \Comment{Solve relaxed RMP and obtain dual variables}
    \State $\mathtt{improved} \leftarrow \texttt{False}$
    \For{each $(r, e) \in \mathcal{B}$}
        \If{$(r, e)$ is dedicated:}
            \For{each $\bar{j}\in\mathcal{J}$}
                \State $(\mathtt{col},\, \bar{c}) \leftarrow \textsc{SolvePRICING}(e, r, \bar{j},\, \boldsymbol{\mu})$ \Comment{MILP or RCSPP formulation}
                \If{$\bar{c} > 0$}
                    \State $\textsc{AddColumn}(\mathtt{col})$;\ \ $\mathtt{improved} \leftarrow \texttt{True}$
                    \State \textbf{continue}
                \EndIf
             \EndFor
         \ElsIf{$(r, e)$ is shared :}
            \State $(\mathtt{col},\, \bar{c}) \leftarrow \textsc{SolvePRICING}(e, r, j,\, \boldsymbol{\mu})$ \Comment{MILP or RCSPP formulation}
            \If{$\bar{c} > 0$}
                \State $\textsc{AddColumn}(\mathtt{col})$;\ \ $\mathtt{improved} \leftarrow \texttt{True}$
                \State \textbf{continue}
            \EndIf
        \EndIf
    \EndFor

\State \Return $\textsc{SolveIntegerRMP}()$ \Comment{Solve integer RMP}

\end{algorithmic}
\end{algorithm}
\section{Computational Analysis}
\label{section:analysis}

All experiments were conducted on a machine equipped with a 96-core Intel(R) Xeon(R) Gold 6338 CPU @ 2.00GHz, with 189GB RAM, under Ubuntu 22.04.
We use \texttt{Python 3.10}, \texttt{Gurobi 12.0.0} for the MILP model, and the \texttt{Pathwyse library} 1.0 \cite{Salani} for the RCSPP formulation. 


\subsection{Instance Generation}
We generate a set of instances designed to evaluate the proposed approach. 
For every instance set we consider a fixed number of specialties, patients, surgeons, and ORs.
In line with the observations of Razali et al.~\cite{Razali}, we adopt a 7-day planning horizon in order to allow the schedule to respond to changes in the surgical waiting list and patient priorities.
To generate the instances, we use the data on surgical procedure distributions fitted in Bernardelli et al.\ \cite{Bernardelli}, which are derived from data of a real OT in Norway \cite{Mannino} and are briefly summarized in Table \ref{tab:data}.
%
\begin{table}[tbp]
    \centering
    \caption{Distribution of surgical specialties, procedure frequencies, and mean durations.
             The column $p$ denotes the proportion of patients per specialty;
             $f_k$ is the relative frequency of patients assigned to procedure $k$
             within the specialty; $\eta_k$ is the mean duration (in minutes) of procedure $k$.}
    \label{tab:data}
    \renewcommand{\arraystretch}{1.3}
    \resizebox{\textwidth}{!}{%
    \begin{tabular}{l l l l}
        \toprule
        \textbf{Specialty} & \textbf{Prop.} $p$ & \textbf{Procedure frequencies} $f_k$ & \textbf{Mean durations} $\eta_k$ (min) \\
        \midrule
        Cardiology       & 0.1432 & $(0.047,\ 0.204,\ 0.172,\ 0.070,\ 0.272,\ 0.083,\ 0.151)$ & $(32.3,\ 53.4,\ 67.2,\ 85.9,\ 107.9,\ 136.3,\ 162.9)$ \\
        Gastroenterology & 0.1876 & $(0.016,\ 0.089,\ 0.224,\ 0.213,\ 0.243,\ 0.182,\ 0.033)$ & $(20.0,\ 41.6,\ 72.6,\ 108.7,\ 159.7,\ 234.0,\ 330.4)$ \\
        Gynecology       & 0.3043 & $(0.206,\ 0.150,\ 0.271,\ 0.328,\ 0.045)$                 & $(24.6,\ 43.8,\ 68.2,\ 124.2,\ 216.4)$ \\
        Orthopedics      & 0.1590 & $(0.040,\ 0.126,\ 0.239,\ 0.535,\ 0.060)$                 & $(32.9,\ 63.9,\ 108.3,\ 174.3,\ 243.6)$ \\
        Urology          & 0.2059 & $(0.158,\ 0.451,\ 0.360,\ 0.031)$                         & $(32.1,\ 58.1,\ 94.7,\ 208.8)$ \\
        \bottomrule
    \end{tabular}}
\end{table}

\subsection{Instance Generation}

We generate a set of benchmark instances to evaluate the proposed approach. Each instance is characterized by a fixed number of specialties, patients, surgeons, and ORs. Following Razali et al.~\cite{Razali}, we adopt a 7-day planning horizon, allowing the schedule to adapt to changes in the surgical waiting list and patient priorities.


We consider two operational settings.
\paragraph{High-Variability LoS ($S_1$)}
This setting represents departments with long and highly variable hospital stays. Weekend staffing is assumed to be reduced by 50\%, so that only half of the available beds can be used in the baseline scenario. Patient LoS follows a distribution with mean $7.1$ days, consistent with~\cite{sdo2023}, and standard deviation $10$ days to represent high variability.

\paragraph{Week Surgery ($S_2$)}
This setting represents departments aiming to discharge all patients before the weekend. Patient pathways are assumed to be simple and standardized, allowing accurate prediction of LoS. Accordingly, LoS is modeled as a discrete uniform random variable over $\{0,\ldots,4\}$ days.

Finally, we consider two OT sizes.
\paragraph{Small ($I_1$)}
Instances with 150 patients, 5 specialties, 30 surgeons, 5 ORs, and 25 beds, representing a small-to-medium operating theater. 

\paragraph{Large ($I_2$)}
Instances with 300 patients, 10 specialties, 60 surgeons, 10 ORs, and 50 beds, representing a medium-to-large operating theater.

Each instance is generated as follows. 
%
Patients are assigned to specialties according to the proportions reported in Table~\ref{tab:data}. Since data are available for only five specialties, the large instances ($I_2$) are generated by duplicating the original specialties. Each patient is assigned a surgical procedure according to the corresponding procedure frequencies in Table~\ref{tab:data}. Surgery durations are set equal to the mean duration of the selected procedure and rounded up to the nearest multiple of 5 minutes.

Patient LoS values are sampled from the distributions described above. Sampled values are rounded up to the nearest integer and truncated at 30 days.

Patient priority scores are generated by drawing an urgency coefficient uniformly from $[0.1,1]$ and scaling it according to the expected surgery duration. Multipliers of 1.15, 1.35, and 1.70 are applied to procedures lasting 60--120, 120--180, and more than 180 minutes, respectively, to compensate for the scheduling disadvantage of longer procedures.

Surgeons are assigned to specialties proportionally to the corresponding waiting-list sizes, ensuring that each specialty has at least one surgeon. Patients are then assigned uniformly at random to surgeons within their specialty.

Turnover times are sampled uniformly between 15 and 40 minutes for consecutive surgeries performed by the same surgeon, and between 40 and 60 minutes otherwise. Values are rounded up to the nearest multiple of 5 minutes, and symmetric turnover times are assumed.

Bed capacities are assigned proportionally to the waiting-list size of each specialty, with a minimum of two beds per specialty.

Each OR operates either as two 4-hour shifts or one 8-hour shift, following \cite{Kianfar}. Half of the ORs adopt each configuration. OR compatibility is generated so that specialties with larger waiting lists are compatible with a larger number of ORs. Under the modified block scheduling policy, 20\% of the ORs are designated as shared, equally distributed between the two shift configurations.

For setting $S_1$, the initial bed occupancy $\hat{\beta}$ is obtained by solving a similar instance with $\Theta_1=\Theta_2=0$ and no initially occupied beds. For setting $S_2$, we assume $\hat{\beta}=0$, since patients are discharged before the weekend.

For each combination of operational setting ($S_1$ or $S_2$) and instance size ($I_1$ or $I_2$), we generate 10 instances.
For each instance, we also consider OR–specialty compatibility and, for setting $S_1$, a number of beds already occupied by patients from the previous week. The complete set of instances is publicly available in the  \href{https://doi.org/10.5281/zenodo.21885087}{online repository}, while an example of an instance is provided in the Supplementary Material.

\subsection{Comparison of the performance of the methods}
\label{section:analysis:performance}

We evaluate the proposed approach by varying the instance size, the scheduling policy, and the parameters controlling weekend bed occupancy ($\Theta_{1_{jd}}\equiv\Theta_1$) and minimum OR utilization ($\Theta_2$). The experiments are organized into four computational analyses comprising 14 tests. Each test is performed on the 10 instances generated for both settings ($S_1$ and $S_2$), resulting in a total of 280 runs.

The four analyses are summarized below.

\begin{description}

\item[Algorithmic comparison.]
We compare the proposed CG approach with the compact MILP formulation on the instance setS $I_1$ AND $I_1$, considering all three scheduling policies and setting $\Theta_1=\Theta_2=0$.

\item[Scalability analysis.]
We evaluate the scalability of the CG approach by comparing its performance on the two instance sizes ($I_1$ and $I_2$) under open scheduling with $\Theta_1=\Theta_2=0$.

\item[Trade-off analysis.]
We investigate the impact of weekend bed availability and minimum OR utilization. Experiments are performed on instance $I_1$ under block scheduling by considering $\Theta_1\in\{0\%,50\%,100\%\}$. The OR utilization threshold is set to $\Theta_2\in\{0\%,65\%\}$ for setting $S_1$ and $\Theta_2\in\{0\%,90\%\}$ for setting $S_2$. The values of $\Theta_2$ are selected to compare the unconstrained case with a representative minimum utilization threshold that remains feasible for all instances.

\item[Policy comparison.]
We compare block, open, and modified block scheduling under both unconstrained ($\Theta_1=\Theta_2=0$) and constrained operating conditions. For the latter, we set $\Theta_1=50\%$, $\Theta_2=65\%$ in $S_1$, and $\Theta_1=50\%$, $\Theta_2=90\%$ in $S_2$.
\end{description}

The results are evaluated using the following performance indicators: OR utilization computed with minimum turnover times (\textbf{OR}${\tilde{\tau}}$), actual OR utilization (\textbf{OR}${\mathrm{real}}$), the number of scheduled high-priority (\textbf{High}), medium-priority (\textbf{Medium}), and low-priority (\textbf{Low}) patients, bed utilization (\textbf{Bed}), objective value (\textbf{Obj}), MIP optimality gap (\textbf{Gap}), and computational time (\textbf{Time}). For the CG approach, computational time is reported separately for CG and the solution of the integer RMP.

The pricing problem is always solved through its RCSPP formulation, as the corresponding MILP formulation proved computationally intractable.
To limit the computational effort, we impose a time limit of 180 minutes for the block scheduling policy, 270 minutes for the modified block scheduling policy, and 360 minutes for the open scheduling policy.  
For the  CG approach, the available time is divided between the  CG phase, where a maximum number of iterations equal to $500$ is also enforced, and the solution of the final RMP.


\subsubsection{Algorithmic comparison}
\label{section:analysis:performance:algorithmic}

In this section, we compare the proposed CG approach with the compact MILP formulation on instances in $I_1$. Results for $I_2$ are omitted because the MILP formulation is unable to find a feasible solution within the time limit, whereas the CG approach consistently produces feasible solutions.
Example solutions for the three scheduling policies are provided in the Supplementary Material.
Table~\ref{tab:milp_small_variable} reports the results for setting $S_1$. Under this setting, the problem is relatively easy to solve, and both approaches reach an optimal (or near-optimal) solution within a short computational time. The CG approach is faster than the MILP formulation under block and modified scheduling, while the two methods exhibit comparable computational times under open scheduling. This is expected, since under open scheduling the pricing problem cannot be decomposed by specialty, increasing the computational effort.
\begin{table}[tbp]
    \centering
    \small
    \caption{Comparison of MILP and CG approaches for instance set $I_1$ and setting $S_1$, with $\Theta_1 = \Theta_2 = 0$.}
    \label{tab:milp_small_variable}
    \renewcommand{\arraystretch}{1.2}
    \resizebox{\textwidth}{!}{%
    \begin{tabular}{l l S[table-format=2.2] S[table-format=2.2] 
                        r r r 
                        S[table-format=2.2] 
                        S[table-format=1.2,table-space-text-post=\%]
                        r r}
        \toprule
        \textbf{Policy} & \textbf{Method} &
        \textbf{OR\,$\tilde{\tau}$\,(\%)} & \textbf{OR real\,(\%)} &
        \textbf{High} & \textbf{Med.} & \textbf{Low} &
        \textbf{Bed\,(\%)} & \textbf{Obj} &
        \textbf{Gap\,(\%)} & \textbf{Time\,(s)} \\
        \midrule
        \multirow{2}{*}{Block} & MILP & 58.46 & 64.85 & 19.70 & 28.40 & 10.90 & 99.46 & 49.25 & 0.00 & 156.16 \\
                       & CG   & 55.71 & 60.64 & 19.60 & 28.40 & 9.20 & 99.46 & 48.56 & 0.13 & 11.24+4.80 \\
\midrule
\multirow{2}{*}{Open} & MILP & 58.62 & 66.16 & 19.70 & 28.40 & 10.90 & 99.52 & 49.25 & 0.00 & 121.58 \\
                      & CG   & 56.68 & 61.83 & 19.80 & 27.90 & 10.70 & 99.18 & 48.97 & 0.04 & 109.71+11.12 \\
\midrule
\multirow{2}{*}{Modified} & MILP & 58.52 & 65.47 & 19.70 & 28.40 & 10.90 & 99.46 & 49.25 & 0.00 & 185.61 \\
                          & CG   & 56.60 & 61.48 & 19.60 & 28.50 & 9.80 & 99.85 & 48.88 & 0.00 & 22.73+5.12 \\
\bottomrule
    \end{tabular}}
\end{table}

The MILP formulation consistently achieves a slightly higher objective value than the CG approach, although the difference is very small. This indicates that the columns generated during the solution of the relaxed RMP are not always sufficient to recover the optimal integer solution.

The MILP formulation attains the same objective value under all three scheduling policies, indicating that the scheduling policy has no impact on the optimal solution for this setting. In contrast, the CG approach shows small differences across policies, with open scheduling performing slightly better than modified and block scheduling. This suggests that the pricing problem generates columns of higher quality under the open scheduling policy.

Finally, OR utilization ranges from $60.64\%$ (CG with block scheduling) to $66.16\%$ (MILP with open scheduling), whereas bed utilization is consistently close to its maximum. This confirms that bed availability, rather than OR capacity, is the primary bottleneck in this operational setting.

Table~\ref{tab:milp_small_week} reports the results for instance $I_1$ under setting $S_2$. In this setting, the problem is considerably more challenging, and the time limit is reached for all MILP instances. Nevertheless, the CG approach consistently achieves higher objective values while requiring less computational time, demonstrating its superior effectiveness on harder instances.
\begin{table}[tbp]
    \centering
    \small
    \caption{Comparison of MILP and CG approaches for  instance $I_1$ and setting $S_2$, with $\Theta_1 = \Theta_2 = 0$.}
    \label{tab:milp_small_week}
    \renewcommand{\arraystretch}{1.2}
    \resizebox{\textwidth}{!}{%
    \begin{tabular}{l l S[table-format=2.2] S[table-format=2.2] 
                        r r r 
                        S[table-format=2.2] 
                        S[table-format=1.2,table-space-text-post=\%]
                        r r}
        \toprule
        \textbf{Policy} & \textbf{Method} &
        \textbf{OR\,$\tilde{\tau}$\,(\%)} & \textbf{OR real\,(\%)} &
        \textbf{High} & \textbf{Med.} & \textbf{Low} &
        \textbf{Bed\,(\%)} & \textbf{Obj} &
        \textbf{Gap\,(\%)} & \textbf{Time\,(s)} \\
        \midrule
         \multirow{2}{*}{Block} & MILP & 81.49 & 94.09 & 27.10 & 44.90 & 14.30 & 95.84 & 71.93 & 3.48 & 10802.98 \\
                           & CG   & 83.22 & 92.88 & 26.90 & 45.20 & 15.60 & 96.32 & 72.56 & 1.34 & 1002.43+12.76 \\
        \midrule
        \multirow{2}{*}{Open} & MILP & 81.71 & 95.44 & 27.30 & 44.60 & 14.70 & 95.28 & 72.12 & 3.21 & 21603.02 \\
                           & CG   & 84.60 & 94.05 & 27.50 & 45.08 & 14.90 & 97.04 & 73.11 & 1.39 & 2789.19+22.35 \\
        \midrule
        \multirow{2}{*}{Modified} & MILP & 80.19 & 93.00 & 26.80 & 44.80 & 14.00 & 94.56 & 71.58 & 3.98 & 
        16203.44 \\
                           & CG & 83.99 & 93.69 & 27.30 & 45.70 & 15.20 & 97.68 & 73.00 & 1.06 & 2190.77+17.43 \\
        \bottomrule
    \end{tabular}}
\end{table}

Both OR and bed utilization are high, indicating a good balance between the two resources. OR utilization reaches $95.44\%$ under the MILP formulation with open scheduling, while bed utilization reaches $97.68\%$ under the CG approach with modified scheduling.

Unlike setting $S_1$, the scheduling policy has a noticeable impact on the solution quality. For both the CG approach and the MILP formulation, open scheduling yields the highest objective values, followed by modified and block scheduling. This behavior is expected, as open scheduling is the most flexible policy, modified scheduling provides intermediate flexibility, and block scheduling is the most restrictive. The additional flexibility allows the model to accommodate more patients within the available OR and bed capacities, resulting in higher objective values.

For the MILP formulation, the block policy slightly outperforms the modified policy. However, this difference is likely due to the larger optimality gap, and both policies could converge to similar objective values if solved to optimality.

\subsubsection{Scalability analysis}
\label{section:analysis:performance:scalability}

In this section, we assess the scalability of the proposed CG approach by comparing its performance on the two instance sizes under the open scheduling policy with $\Theta_1=\Theta_2=0$.

Table~\ref{tab:size} reports the results. As expected, increasing the instance size increases both the computational time and the optimality gap. Nevertheless, the proposed approach remains capable of solving the larger instances within acceptable computational times while maintaining relatively small optimality gaps. The objective value and the number of scheduled patients approximately double from $I_1$ to $I_2$, reflecting the doubled problem size.
\begin{table}[tbp]
    \centering
    \small
    \caption{Instance size comparison with $\Theta_1 = 0\%$ and $\Theta_2 = 0\%$ under the open scheduling policy.}
    \label{tab:size}
    \renewcommand{\arraystretch}{1.2}
    \resizebox{\textwidth}{!}{%
    \begin{tabular}{l l
                    S[table-format=2.2,table-space-text-post=\%]
                    S[table-format=2.2,table-space-text-post=\%]
                    r r r
                    S[table-format=2.2,table-space-text-post=\%]
                    S[table-format=3.2]
                    S[table-format=1.2,table-space-text-post=\%]
                    r}
        \toprule
         & &
        \textbf{OR\,$\tilde{\tau}$\,(\%)} & \textbf{OR real\,(\%)} &
        \textbf{High} & \textbf{Med.} & \textbf{Low} &
        \textbf{Bed\,(\%)} & \textbf{Obj} &
        \textbf{Gap\,(\%)} & \textbf{Time\,(s)} \\
        \midrule
       $I_1$ & $S_1$ & 56.68 & 61.83 & 19.80 & 27.90 & 10.70 & 99.18 & 48.97 & 0.04 & 109.71+11.12 \\
$I_2$ & $S_1$ & 54.83 & 59.07 & 38.40 & 56.20 & 20.20 & 98.78 & 96.12 & 1.23 & 1301.20+33.99 \\
$I_1$ & $S_2$ & 84.60 & 94.05 & 27.50 & 45.08 & 14.90 & 97.04 & 73.11 & 1.39 & 2789.19+22.35 \\
$I_2$ & $S_2$ & 82.05 & 90.43 & 50.40 & 89.80 & 31.00 & 94.72 & 140.06 & 3.63 & 8282.86+37.46 \\
        \bottomrule
    \end{tabular}}
\end{table}

The computational difficulty is strongly influenced by the operational setting. Setting $S_1$ is consistently easier to solve than $S_2$, achieving smaller optimality gaps in considerably shorter computational times. This difference is explained by the resource utilization patterns. In setting $S_1$, OR utilization remains relatively low (up to $61.83\%$), whereas bed utilization is almost saturated (up to $99.18\%$), indicating that bed availability is the primary bottleneck. In contrast, setting $S_2$ exhibits a more balanced use of ORs and beds, with OR utilization reaching $94.05\%$ and bed utilization remaining high (up to $97.04\%$). As a result, both resources become binding, making the scheduling problem harder.

\subsubsection{Trade-off analysis}
\label{section:analysis:performance:resource}
For this set of experiments, we investigate the trade-off between the three objectives of the optimization model. By varying $\Theta_1$ and $\Theta_2$ in the single-objective MILP formulation (Section~\ref{subsec:math}), different Pareto-optimal solutions can be generated, allowing us to evaluate the interaction between patient priority, OR utilization, and weekend bed occupancy.

Tables~\ref{tab:tradeoffs_variable} and~\ref{tab:tradeoffs_week} report the results for settings $S_1$ and $S_2$, respectively. Italicized rows correspond to averages computed over the subset of feasible instances, whereas the second part of each table reports averages restricted to the instances that are feasible for all parameter combinations. 
\begin{table}[tbp]
    \centering
    \small
    \caption{Results of the trade-off analysis under the variable LoS setting $S_1$ using the block scheduling policy. For $\Theta_1 = 0$, 5 of the 10 instances were feasible, whereas for $\Theta_1 = 50$, 9 instances were feasible. The upper part of the table reports results for all feasible instances under each parameter combination. The lower part reports results only for the subset of instances that are feasible across all parameter combinations, ensuring a consistent basis for comparison.}
    \label{tab:tradeoffs_variable}
    \renewcommand{\arraystretch}{1.2}
    \resizebox{\textwidth}{!}{%
    \begin{tabular}{S[table-format=3.0] S[table-format=2.0]
                    S[table-format=2.2,table-space-text-post=\%]
                    S[table-format=2.2,table-space-text-post=\%]
                    r r r
                    S[table-format=3.2,table-space-text-post=\%]
                    S[table-format=3.2]
                    S[table-format=1.2,table-space-text-post=\%]
                    r}
        \toprule
        \textbf{$\Theta_1$} & \textbf{$\Theta_2$} &
        \textbf{OR\,$\tilde{\tau}$\,(\%)} & \textbf{OR real\,(\%)} &
        \textbf{High} & \textbf{Med.} & \textbf{Low} &
        \textbf{Bed\,(\%)} & \textbf{Obj} &
        \textbf{Gap\,(\%)} & \textbf{Time\,(s)} \\
        \midrule
        0 & 0 &  55.71 & 60.64 & 19.60 & 28.40 & 9.20 & 99.46 & 48.56 & 0.13 & 11.24+4.80 \\
\textit{0} & \textit{65} & \textit{65.53} & \textit{72.28} & \textit{20.40} & \textit{28.40} & \textit{13.60} & \textit{99.32} & \textit{50.46} & \textit{0.45} & \textit{703.78+8.99} \\
\midrule
50 & 0 & 60.33 & 66.07 & 23.20 & 29.00 & 8.90 & 98.87 & 53.64 & 0.42 & 134.08+5.96 \\
50 & 65 & \textit{66.41} & \textit{73.03} & \textit{24.00} & \textit{28.22} & \textit{11.89} & \textit{99.5}1 & \textit{54.63} & \textit{0.30} & \textit{239.07+11.05} \\
\midrule
100 & 0 & 63.94 & 70.38 & 24.70 & 30.30 & 9.30 & 98.23 & 56.71 & 0.24 & 20.48+5.72 \\
100 & 65 & 66.49 & 73.72 & 25.10 & 29.60 & 11.30 & 98.76 & 56.91 & 0.25 & 290.73+9.20 \\
\midrule
\midrule
        0 & 0 &  58.75 & 64.38 & 21.20 & 29.00 & 9.20 & 98.78 & 50.90 & 0.26 & 11.78 + 4.99 \\
0 & 65 & 65.53 & 72.28 & 20.40 & 28.40 & 13.60 & 99.32 & 50.46 & 0.45 & 703.78+8.99 \\
\midrule
50 & 0 & 63.19 & 69.47 & 20.40 & 28.40 & 13.60 & 99.32 & 55.91 & 0.60 & 254.17+5.98 \\
50 & 65 & 67.46 & 74.40 & 24.60 & 29.20 & 12.00 & 99.62 & 55.78 & 0.39 & 41.51+11.18 \\
\midrule
100 & 0 & 65.86 & 72.96 & 26.00 & 31.00 & 8.60 & 97.89 & 58.46 & 0.24 & 23.39+7.57 \\
100 & 65 & 67.30 & 74.72 & 26.40 & 30.60 & 10.60 & 98.83 & 58.86 & 0.11 & 28.60+8.76 \\

        \bottomrule
    \end{tabular}}
\end{table}

For setting $S_1$, OR $\tilde{\tau}$ without minimum utilization requirements ranges from $55.71\%$ for $\Theta_1=0$ to $63.94\%$ for $\Theta_1=100$. Increasing $\Theta_2$ makes the problem substantially harder: with $\Theta_1=0$, half of the instances become infeasible, whereas only one instance is infeasible for $\Theta_1=50$.

Considering only the instances that remain feasible for all parameter combinations, the objective value decreases as $\Theta_2$ increases when $\Theta_1\in\{0,50\}$, showing that the utilization requirement forces the model to schedule lower-priority patients. For $\Theta_1=100$, however, a slight increase in the objective value is observed. Although counterintuitive, this suggests that the additional utilization constraint may generate columns that are more effective for the integer RMP. Overall, minimum utilization requirements improve CG efficiency when sufficient resources are available, but may reduce solution quality and increase infeasibility under tighter resource conditions.

\begin{table}[tb]
    \centering
    \small
    \caption{Comparison of resource trade-offs under the variable LoS setting $S_2$ using the block scheduling policy. For $\Theta_1 = 0$, 2 of the 10 instances were feasible, whereas for $\Theta_1 = 50$, 9 instances were feasible. The upper part of the table reports results for all feasible instances under each parameter combination. The lower part reports results only for the subset of instances that are feasible across all parameter combinations, ensuring a consistent basis for comparison.}
    \label{tab:tradeoffs_week}
    \renewcommand{\arraystretch}{1.2}
    \resizebox{\textwidth}{!}{%
    \begin{tabular}{S[table-format=3.0] S[table-format=2.0]
                    S[table-format=2.2,table-space-text-post=\%]
                    S[table-format=2.2,table-space-text-post=\%]
                    r r r
                    S[table-format=3.2,table-space-text-post=\%]
                    S[table-format=3.2]
                    S[table-format=1.2,table-space-text-post=\%]
                    r}
        \toprule
        \textbf{$\Theta_1$} & \textbf{$\Theta_2$} &
        \textbf{OR\,$\tilde{\tau}$\,(\%)} & \textbf{OR real\,(\%)} &
        \textbf{High} & \textbf{Med.} & \textbf{Low} &
        \textbf{Bed\,(\%)} & \textbf{Obj} &
        \textbf{Gap\,(\%)} & \textbf{Time\,(s)} \\
          \midrule
        0 & 0 & 83.22 & 92.88 & 26.90 & 45.20 & 15.60 & 96.32 & 72.56 & 1.34 & 1002.43+12.76  \\
\textit{0} & \textit{90}&  \textit{90.06}  & \textit{97.73}  &\textit{ 29.50} & \textit{45.50} & \textit{15.50} & \textit{92.00} & \textit{75.34} & \textit{2.79} & \textit{1022.13+12.64} \\
\midrule
50 & 0 & 85.72 & 94.83 & 28.30 & 50.50 & 12.90 & 94.08 & 77.43 & 2.28 & 1481.54+12.76 \\
\textit{50} & \textit{90} & \textit{90.25} & \textit{98.39} & \textit{27.22} & \textit{48.89} & \textit{16.44} & \textit{92.38} & \textit{75.97} & \textit{3.39} & \textit{1117.91+12.7}2 \\
\midrule
100 & 0 & 87.24 & 96.68 & 28.60 & 53.90 & 13.20 & 92.00 & 80.59 & 1.35 & 1446.87+12.65 \\
100 & 90 & 90.24 & 98.35 & 28.40 & 51.80 & 15.50 & 89.43 & 79.18 & 2.14 & 1314.29+12.74 \\
\midrule
\midrule
0 & 0& 86.98 & 97.12& 31.00 & 46.00 & 13.00 & 96.40 & 77.42 & 0.93 & 1023.56+12.80 \\  
0 & 90 & 90.06  & 97.73 & 29.50 & 45.50 & 15.50 & 92.00 & 75.34 & 2.79 & 1022.13+12.64 \\ 
\midrule
50 & 0 & 88.40 & 96.62 & 32.00 & 48.50 & 12.00 & 92.95 & 80.35 & 2.32 & 1089.58+12.27 \\ 
50 & 90& 90.10 & 97.33 & 31.50 & 48.00 & 13.50 & 90.60 & 80.01 & 2.25 & 106.91+12.43 \\ 
\midrule
100 & 0 & 88.48 & 97.42 & 32.50 & 52.50 & 9.50 & 90.29 & 83.07 & 0.93 & 3049.54+13.22 \\ 
100 & 90& 90.23 & 98.23 & 33.00 & 49.50 & 12.00 & 90.00 & 82.78 & 0.84 & 1185.13+13.09 \\
        \bottomrule
    \end{tabular}}
\end{table}

Setting $S_2$ exhibits similar behavior. Without utilization requirements, OR $\tilde{\tau}$ ranges from $83.22\%$ to $87.24\%$. Imposing $\Theta_2=90\%$ renders eight instances infeasible for $\Theta_1=0$ and one instance infeasible for $\Theta_1=50$. For the instances that remain feasible, the objective value consistently decreases as $\Theta_2$ increases, confirming that higher utilization requirements force the model to prioritize OR efficiency over patient priority. In return, actual OR utilization increases, reaching up to $98.23\%$.

\subsubsection{Policy comparison}
\label{section:analysis:performance:policy}
In this section, we compare the three scheduling policies under different operating conditions to evaluate their impact on solution quality, resource utilization, and computational performance.

\begin{table}[tbp]
    \centering
    \small
    \caption{Policy comparison with $\Theta_1 = 0\%$ and $\Theta_2 = 0\%$ for instance $I_2$ and setting $S_1$.}
    \label{tab:policy_00_variable}
    \renewcommand{\arraystretch}{1.2}
    \resizebox{\textwidth}{!}{%
    \begin{tabular}{l
                    S[table-format=2.2,table-space-text-post=\%]
                    S[table-format=2.2,table-space-text-post=\%]
                    r r r
                    S[table-format=2.2,table-space-text-post=\%]
                    S[table-format=3.2]
                    S[table-format=1.2,table-space-text-post=\%]
                    r}
        \toprule
        \textbf{Policy} &
        \textbf{OR\,$\tilde{\tau}$\,(\%)} & \textbf{OR real\,(\%)} &
        \textbf{High} & \textbf{Med.} & \textbf{Low} &
        \textbf{Bed\,(\%)} & \textbf{Obj} &
        \textbf{Gap\,(\%)} & \textbf{Time\,(s)} \\
        \midrule
        Block & 55.94 & 60.21 & 38.50 & 56.20 & 21.60 & 99.15 & 96.51 & 0.72 & 60.73+23.58 \\
        Open & 54.83 & 59.07 & 38.40 & 56.20 & 20.20 & 98.78 & 96.12 & 1.23 & 1301.20+33.99 \\
        Modified & 56.34 & 60.74 & 38.30 & 56.20 & 23.00 & 99.32 & 96.70 & 0.52 & 134.13+27.53 \\
        \bottomrule
    \end{tabular}}
\end{table}

Table~\ref{tab:policy_00_variable} reports the results for setting $S_1$ with $\Theta_1=\Theta_2=0$. The modified scheduling policy achieves the highest average objective value, followed by the block and open policies. However, the differences are very small, suggesting that the three scheduling policies are nearly equivalent under this setting. The observed differences are therefore likely attributable to the computational performance of the solution approach rather than to the scheduling policy itself.

\begin{table}[tbp]
    \centering
    \small
    \caption{Policy comparison with $\Theta_1 = 50\%$ and $\Theta_2 = 65\%$ for instance $I_2$ and setting $S_1$. For all scheduling policies, one instance results infeasible. The reported results are based on the 9 feasible instances.}
    \label{tab:policy_5095_variable}
    \renewcommand{\arraystretch}{1.2}
    \resizebox{\textwidth}{!}{%
    \begin{tabular}{l
                    S[table-format=2.2,table-space-text-post=\%]
                    S[table-format=2.2,table-space-text-post=\%]
                    r r r
                    S[table-format=2.2,table-space-text-post=\%]
                    S[table-format=3.2]
                    S[table-format=1.2,table-space-text-post=\%]
                    r}
        \toprule
        \textbf{Policy} &
        \textbf{OR\,$\tilde{\tau}$\,(\%)} & \textbf{OR real\,(\%)} &
        \textbf{High} & \textbf{Med.} & \textbf{Low} &
        \textbf{Bed\,(\%)} & \textbf{Obj} &
        \textbf{Gap\,(\%)} & \textbf{Time\,(s)} \\
        \midrule
        Block & 65.50 & 71.74  & 44.78 & 57.22 & 25.78 & 98.85 & 106.04 &1.18& 5253.70+28.51 \\
        Open & 65.20 & 71.13 & 44.78 & 57.00 & 26.11 & 98.61 & 106.00 & 1.33 & 7950.48+36.72 \\
        Modified & 65.67 & 71.91 & 44.89 & 57.11 & 27.00 & 98.96 & 106.24 & 1.11 & 2595.47+31.70 \\
        \bottomrule
    \end{tabular}}
\end{table}

Table~\ref{tab:policy_5095_variable} reports the results for setting $S_1$ with $\Theta_1=50\%$ and $\Theta_2=65\%$. One of the ten instances is infeasible under all three scheduling policies. Although the overall performance of the three policies is similar, the modified scheduling policy consistently provides the best computational performance, achieving the highest objective value together with the smallest optimality gap.

\begin{table}[tbp]
    \centering
    \small
    \caption{Policy comparison with $\Theta_1 = 0\%$ and $\Theta_2 = 0\%$ for instance $I_2$ and setting $S_2$.}
    \label{tab:policy_00_week}
    \renewcommand{\arraystretch}{1.2}
     \resizebox{\textwidth}{!}{%
    \begin{tabular}{l
                    S[table-format=2.2,table-space-text-post=\%]
                    S[table-format=2.2,table-space-text-post=\%]
                    r r r
                    S[table-format=2.2,table-space-text-post=\%]
                    S[table-format=3.2]
                    S[table-format=1.2,table-space-text-post=\%]
                    r}
        \toprule
        \textbf{Policy} &
        \textbf{OR\,$\tilde{\tau}$\,(\%)} & \textbf{OR real\,(\%)} &
        \textbf{High} & \textbf{Med.} & \textbf{Low} &
        \textbf{Bed\,(\%)} & \textbf{Obj} &
        \textbf{Gap\,(\%)} & \textbf{Time\,(s)} \\
        \midrule
       Block & 81.55 & 90.27 & 50.40 & 90.60 & 31.10 & 95.64 & 140.59 & 2.46 & 3418.56+29.63 \\
Open & 82.05 & 90.43 & 50.40  & 89.80 & 31.00 & 94.72 & 140.06 &  3.63  & 8282.86+37.46 \\
Modified & 82.36 & 90.94 & 50.50 & 90.90 & 31.60 & 95.72 & 141.07 & 2.53 & 5636.52+33.30 \\
        \bottomrule
    \end{tabular}}
\end{table}

Table~\ref{tab:policy_00_week} reports the results for setting $S_2$ with $\Theta_1=\Theta_2=0$. The modified scheduling policy achieves the highest objective value, followed by the block policy, whereas the open policy performs worst despite the additional computational time allocated to this policy. This behavior is consistent with the higher computational complexity of the open scheduling formulation. OR and bed utilization follow the same ranking across the three scheduling policies.

\begin{table}[tbp]
    \centering
    \small
    \caption{Policy comparison for instance $I_2$ under setting $S_2$, with $\Theta_1 = 50\%$ and $\Theta_2 = 90\%$. For the block and modified scheduling policies, 8 of the 10 instances were feasible, whereas for the open scheduling policy, 5 instance were fesible. The results are reported as the mean over the five instances that are feasible under all three policies.}
    \label{tab:policy_5095_week}
    \renewcommand{\arraystretch}{1.2}
     \resizebox{\textwidth}{!}{%
    \begin{tabular}{l
                    S[table-format=2.2,table-space-text-post=\%]
                    S[table-format=2.2,table-space-text-post=\%]
                    r r r
                    S[table-format=2.2,table-space-text-post=\%]
                    S[table-format=3.2]
                    S[table-format=1.2,table-space-text-post=\%]
                    r}
        \toprule
        \textbf{Policy} &
        \textbf{OR\,$\tilde{\tau}$\,(\%)} & \textbf{OR real\,(\%)} &
        \textbf{High} & \textbf{Med.} & \textbf{Low} &
        \textbf{Bed\,(\%)} & \textbf{Obj} &
        \textbf{Gap\,(\%)} & \textbf{Time\,(s)} \\
        \midrule
     Block & 90.03 & 97.51 & 55.00 & 94.20 & 35.60 & 92.23 & 150.07 & 4.93 & 4922.70+31.41 \\
Open & 90.09 & 97.93 & 52.60 & 96.20 & 31.60 & 88.51 &  147.94 & 9.26 & 15069.08+39.49  \\
Modified & 90.03 & 97.97 & 54.20 & 96.60 & 34.40 & 90.47 & 150.75 & 6.10 & 9820.95+35.24  \\
        \bottomrule
    \end{tabular}}
\end{table}

Table~\ref{tab:policy_5095_week} reports the results for setting $S_2$ with $\Theta_1=50\%$ and $\Theta_2=90\%$. Considering only the instances that are feasible under all three scheduling policies, the modified scheduling policy consistently achieves the highest objective value within the imposed time and iteration limits.

An interesting observation is that the infeasible instances differ between the block and modified scheduling policies. In both cases, infeasibility occurs because the integer RMP fails to identify a solution with zero slack. This indicates that the two policies fail on different subsets of instances, suggesting that they exhibit complementary feasibility characteristics.


\subsection{Performance Analysis under Uncertainty}
\label{section:analysis:uncertainty}
In this section, we evaluate the robustness of the proposed approach under uncertain surgery durations and patient LoS.

We assume that the surgical schedule cannot be modified during execution. Consequently, surgeries must be performed on their assigned day and cannot start before their planned start times. Delays caused by longer-than-expected procedures therefore result in OR overtime, whereas shorter procedures generate idle time. Similarly, if the realized bed occupancy exceeds the available capacity, patients are assumed to be temporarily accommodated in beds belonging to other specialties.

Patients are scheduled according to the schedule produced by the optimization model, with consecutive surgeries separated only by the corresponding turnover times. Surgery durations are sampled from lognormal distributions fitted by Bernardelli et al.~\cite{Bernardelli}, using the mean values reported in Table~\ref{tab:data} and the standard deviations reported in Table~\ref{tab:uncertainty}.

\begin{table}[tbp]
    \centering
    \caption{Standard deviation $\sigma_k$ for each procedure type.}
    \label{tab:uncertainty}
    \begin{tabular}{l l}
        \toprule
        \textbf{Specialty} &  \textbf{Standard Deviation} $\sigma_k$ (min) \\
        \midrule
        Cardiology    & $(4.9,\ 10.8,\ 13.6,\ 17.4,\ 21.8,\ 20.6,\ 24.6)$ \\
        Gastroenterology & $(4.0,\ 8.4,\ 14.7,\ 16.4,\ 24.1,\ 35.3,\ 49.8)$ \\
        Gynecology   & $(5.0,\ 6.6,\ 10.3,\ 25.1,\ 32.6)$ \\
        Orthopedics  & $(6.7,\ 9.6,\ 21.9,\ 26.3,\ 49.2)$ \\
        Urology   & $(6.5,\ 11.7,\ 19.1,\ 31.5)$ \\
        \bottomrule
    \end{tabular}
\end{table}

To generate realistic LoS values, we sample from a lognormal distribution \cite{Marazzi} with mean equal to the estimated LoS and standard deviation equal to $10\%$ of the mean.

For each schedule, we generate 1,000 scenarios by resampling surgery durations and LoS values. The resulting schedules are evaluated using the following performance indicators: average patient waiting time (\textbf{WT}, min), average OR idle time (\textbf{Idle time}, min), average OR overtime (\textbf{Overtime}, min), actual OR utilization (\textbf{OR real}, \%), average excess bed occupancy (\textbf{Excess}), and average bed utilization (\textbf{Bed}, \%).
The experimental settings are the same as those considered in Sections~\ref{section:analysis:performance:resource} and~\ref{section:analysis:performance:policy}, and all reported values are averaged over the 1,000 scenarios.

\subsubsection{Trade-off analysis}

\begin{table}[tbp]
    \centering
    \small
    \caption{Resource trade-off comparison under uncertainty for instance $I_1$ under the block scheduling policy. The reported results represent the average values computed over all instances that are feasible for each parameter combination.}
\label{tab:unc:tradeoffs}
    \renewcommand{\arraystretch}{1.2}
    \resizebox{\textwidth}{!}{%
    \begin{tabular}{l l r r r r r r r r r r r r r r}
        & & & \multicolumn{6}{c}{$S_1$} & & \multicolumn{6}{c}{$S_2$} \\
        \cmidrule(lr){4-9} \cmidrule(lr){11-16}
        \textbf{$\Theta_1$} & \textbf{$\Theta_2$} & &
        \textbf{WT}  &
        \textbf{Idle} & \textbf{Over} & \textbf{OR real} & \textbf{Excess} & \textbf{Beds (\%)}& &
        \textbf{WT} &
        \textbf{Idle} & \textbf{Over} & \textbf{OR real} & \textbf{Excess} & \textbf{Beds (\%)} \\
        \midrule
        0 & 0 &  & 0.33 & 8.46 & 2.19 & 0.66 & 3.82 & 0.98 & & 0.69 & 7.62 & 7.49 &0.99 & 7.62 & 0.96\\
        & 65/90  & & 0.54 & 7.30 & 3.55 & 0.74& 3.47& 0.98 &  & 1.03 & 7.16 & 7.84 & 1.00 & 6.28 & 0.91 \\
        \midrule
        \multirow{2}{*}{50} & 0 & & 0.41 & 6.17 & 1.94 &0.71&3.13&0.97&  & 0.92 & 8.25 & 6.86 &0.99& 5.01 & 0.91 \\
                           & 65/90  & & 0.58 & 7.10 & 3.32 & 0.76 &2.95 & 0.97 &  & 1.02 & 7.01 & 7.59 & 0.99 &4.02 & 0.89\\
        \midrule
        \multirow{2}{*}{100} & 0 &  & 0.39 & 9.28 & 2.68 & 0.74& 2.50 & 0.94 & & 0.77 & 10.09 & 6.88 &0.99 &3.41 & 0.87 \\
                           & 65/90   &  & 0.48 & 9.11 & 2.50 & 0.76 & 2.49 & 0.95 &  & 1.04 & 8.20 & 7.73 & 1.00&  3.17 & 0.86\\
        \bottomrule
    \end{tabular}}
\end{table}

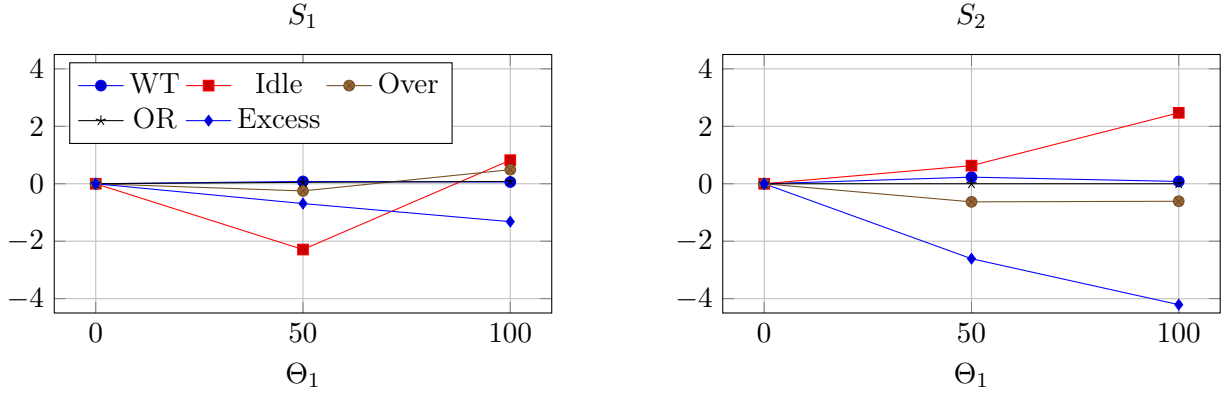
\begin{figure}[h!]
\centering
\caption{Impact of increasing $\Theta_1$ on performance metrics for settings $S_1$ and $S_2$. Values are reported as changes with respect to $\Theta_1=0$, considering $\Theta_2=0$ .}
\label{fig:beds_impact}
\begin{minipage}{0.48\textwidth}
\centering
\begin{tikzpicture}
\begin{axis}[
    width=\textwidth,
    height=5cm,
    xlabel={$\Theta_1$},
    ylabel={},
    xtick={0,50,100},
    legend pos=north west,
    legend columns=3,
    grid=major,
    ymin=-4.5,
    ymax=4.5,
    title={$S_1$}
]

\addplot coordinates {
(0,0)
(50,0.08)
(100,0.06)
};

\addplot coordinates {
(0,0)
(50,-2.29)
(100,0.82)
};

\addplot coordinates {
(0,0)
(50,-0.25)
(100,0.49)
};

\addplot coordinates {
(0,0)
(50,0.05)
(100,0.08)
};

\addplot coordinates {
(0,0)
(50,-0.69)
(100,-1.32)
};

\legend{WT,Idle,Over,OR,Excess}

\end{axis}
\end{tikzpicture}
\end{minipage}
\hfill
\begin{minipage}{0.48\textwidth}
\centering
\begin{tikzpicture}
\begin{axis}[
    width=\textwidth,
    height=5cm,
    xlabel={$\Theta_1$},
    ylabel={},
    xtick={0,50,100},
    grid=major,
    ymin=-4.5,
    ymax=4.5,
    title={$S_2$}
]


\addplot coordinates {
(0,0)
(50,0.23)
(100,0.08)
};

\addplot coordinates {
(0,0)
(50,0.63)
(100,2.47)
};

\addplot coordinates {
(0,0)
(50,-0.63)
(100,-0.61)
};

\addplot coordinates {
(0,0)
(50,0)
(100,0)
};

\addplot coordinates {
(0,0)
(50,-2.61)
(100,-4.21)
};


\end{axis}
\end{tikzpicture}
\end{minipage}
\end{figure}
In this Section, we analyze the impact of increased weekend bed availability and minimum OR occupation requirements on the model performance under uncertainty realization. 
As shown in Figure~\ref{fig:beds_impact} and Table \ref{tab:unc:tradeoffs}, increasing weekend bed availability has a limited effect on waiting times and OR occupation, which remain largely unchanged. Excess bed occupancy decreases across all scenarios, while overtime exhibits a slight reduction only under setting $S_2$. The impact on idle time is less straightforward: under setting $S_1$, idle time first decreases and then increases, whereas under setting $S_2$ it increases consistently as weekend bed availability grows.

\begin{figure}[h!]
\centering
\caption{Impact of the occupancy constraint. The bars represent the difference between the values of the indicators (WT, Idle, Over, and Excess) obtained under a minimum occupancy requirement and those obtained without any occupancy constraint.}
\label{fig:occuancy_impact}
\begin{minipage}{0.48\textwidth}
\centering
\begin{tikzpicture}
\begin{axis}[
    ybar,
    bar width=6pt,
    width=8cm,
    height=5.5cm,
    ylabel={$f(\Theta_2=65/90)-f(\Theta_2=0)$},
    symbolic x coords={0,50,100},
    xtick=data,
    xlabel={$\Theta_1$},
    legend pos=south east,
    legend columns = 3,
    enlarge x limits=0.2,
    ymajorgrids,
    ymin=-2,
    ymax=2,
    title={$S_1$}
]

\addplot coordinates {
(0,0.21)
(50,0.17)
(100,0.09)
};

\addplot coordinates {
(0,-1.16)
(50,0.93)
(100,-0.17)
};

\addplot coordinates {
(0,1.36)
(50,1.38)
(100,-0.18)
};

\addplot coordinates {
(0,0.08)
(50,0.05)
(100,0.02)
};

\addplot coordinates {
(0,-0.35)
(50,-0.18)
(100,-0.01)
};

\legend{WT,Idle,Over,OR,Excess}

\end{axis}
\end{tikzpicture}
\end{minipage}
\hfill
\begin{minipage}{0.48\textwidth}
\centering
\begin{tikzpicture}
\begin{axis}[
    ybar,
    bar width=6pt,
    width=8cm,
    height=5.5cm,
    ylabel={},
    symbolic x coords={0,50,100},
    xtick=data,
    xlabel={$\Theta_1$},
    legend pos=north east,
    enlarge x limits=0.2,
    ymajorgrids,
    ymin=-2,
    ymax=2,
    title={$S_2$}
]

\addplot coordinates {
(0,0.34)
(50,0.1)
(100,0.27)
};

\addplot coordinates {
(0,-0.46)
(50,-1.24)
(100,-1.89)
};

\addplot coordinates {
(0,0.35)
(50,0.73)
(100,0.85)
};

\addplot coordinates {
(0,0.01)
(50,0.0)
(100,0.01)
};

\addplot coordinates {
(0,-1.34)
(50,-0.99)
(100,-0.24)
};


\end{axis}
\end{tikzpicture}
\end{minipage}
\end{figure}
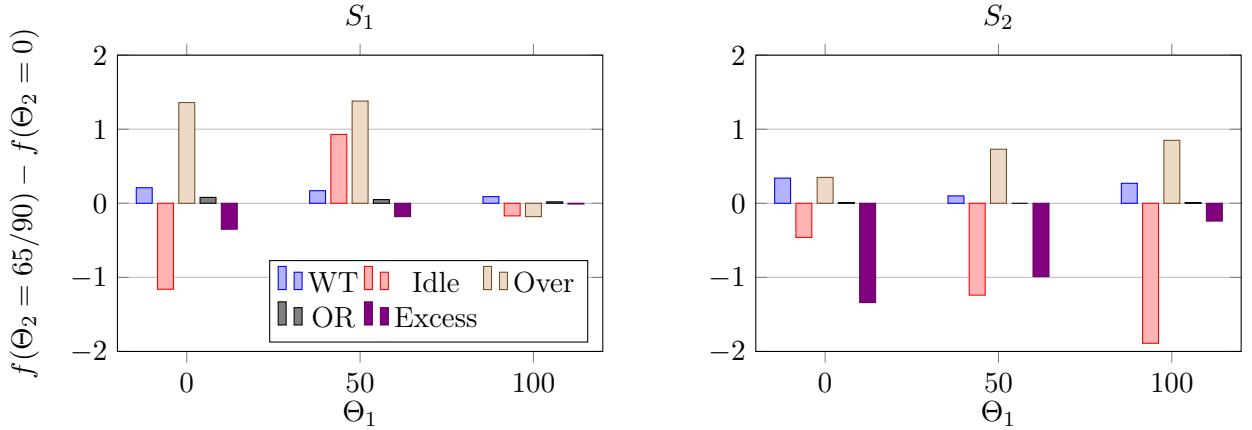

As shown in Figure~\ref{fig:occuancy_impact} and Table \ref{tab:unc:tradeoffs}, the impact of imposing a minimum utilization requirement depends on the level of bed availability. Under setting $S_1$, waiting times increase slightly across all bed configurations. Idle time decreases for $\Theta_1 = 0$ and $\Theta_1 = 100$, but increases for $\Theta_1 = 50$. Overtime increases when $\Theta_1 = 0$ and $\Theta_1 = 50$, whereas a slight reduction is observed for $\Theta_1 = 100$. Finally, excess bed occupancy decreases for all bed configurations; however, this reduction becomes less pronounced as weekend bed availability increases.
In setting $S_2$, the introduction of a minimum utilization requirement results in higher waiting times and overtime. Conversely, idle time and excess bed occupancy decrease, indicating a more intensive use of available resources at the expense of patient waiting times and staff workload.
OR occupation increases slightly but remains approximately stable. This suggests that the constraint on OR occupation loses importance when uncertainty is introduced, and that maximizing patient priority is a sufficiently strong proxy for this objective.
\subsubsection{Policy Comparison}
In this Section, we analyze the impact of the policies on the model performance under uncertainty realization. 
\begin{table}[tbp]
    \small
    \caption{Policy Comparison with $\Theta_1=\Theta_2=0$ under uncertainty for instance $I_2$}
    \label{tab:unc:policy_00}   \renewcommand{\arraystretch}{1.2}
    \resizebox{\textwidth}{!}{%
    \begin{tabular}{l r r r r r r r r r r r r r r}
        & &  \multicolumn{6}{c}{$S_1$} & & \multicolumn{6}{c}{$S_2$} \\
        \cmidrule(lr){3-8} \cmidrule(lr){10-15}
        \textbf{Policy} & &
        \textbf{WT}  &
        \textbf{Idle} & \textbf{Over} & \textbf{OR real} & \textbf{Excess} & \textbf{Beds (\%)}& &
        \textbf{WT} &
        \textbf{Idle}& \textbf{Over} & \textbf{OR real} & \textbf{Excess} & \textbf{Beds (\%)}\\
         \midrule 
        Block && 0.26 & 4.71 & 1.72 & 0.61&8.13 & 0.97 &  & 0.72  & 6.47 & 5.42 & 0.93 &13.84 & 0.95 \\
Open  &&  0.12 & 4.08 & 1.72 & 0.59 &7.77 & 0.97  &  & 0.74 & 6.36 & 6.50 & 0.93 &13.45&0.94 \\
Modified   && 0.23 & 4.52 & 1.68 & 0.61 & 8.00 & 0.97 & & 0.64 & 6.21 & 6.10 & 0.93 & 13.66 & 0.95 \\
        \bottomrule
    \end{tabular}}
\end{table}

\begin{figure}[h!]
\caption{Comparison of uncertainty indicators across the three scheduling policies, considering $\Theta_1 = \Theta_2=0$ for instance $I_2$.}
\label{fig:uncertainty_policy_00}
\centering
\begin{tikzpicture}

\begin{groupplot}[
group style={
    group size=2 by 1,
    horizontal sep=2cm,
},
enlarge x limits=0.15,
ybar,
width=7cm,
height=4.5cm,
symbolic x coords={WT,Idle,Over,OR,Excess},
xtick=data,
ymin=-0.1,
ymax=15,
ymajorgrids=true,
axis x line*=bottom,
axis y line*=left,
legend style={
    at={(1.15,1.05)},
    anchor=south,
    legend columns=3,
},
]

\nextgroupplot[
title={$S_1$},
ybar,
bar width=7pt,
]

\addplot coordinates {
(WT,0.26)
(Idle,4.71)
(Over,1.72)
(OR,0.61)
(Excess,8.13)
};

\addplot coordinates {
(WT,0.12)
(Idle,4.08)
(Over,1.72)
(OR,0.59)
(Excess,7.77)
};

\addplot coordinates {
(WT,0.23)
(Idle,4.52)
(Over,1.68)
(OR,0.61)
(Excess,8.00)
};

\legend{Block,Open,Modified}

\nextgroupplot[
title={$S_2$},
ybar,
bar width=7pt,
]

\addplot coordinates {
(WT,0.72)
(Idle,6.47)
(Over,5.42)
(OR, 0.93)
(Excess,13.84)
};

\addplot coordinates {
(WT,0.74)
(Idle,6.36)
(Over,6.50)
(OR, 0.93)
(Excess,13.45)
};

\addplot coordinates {
(WT,0.64)
(Idle,6.21)
(Over,6.10)
(OR, 0.93)
(Excess,13.66)
};

\end{groupplot}

\end{tikzpicture}
\end{figure}
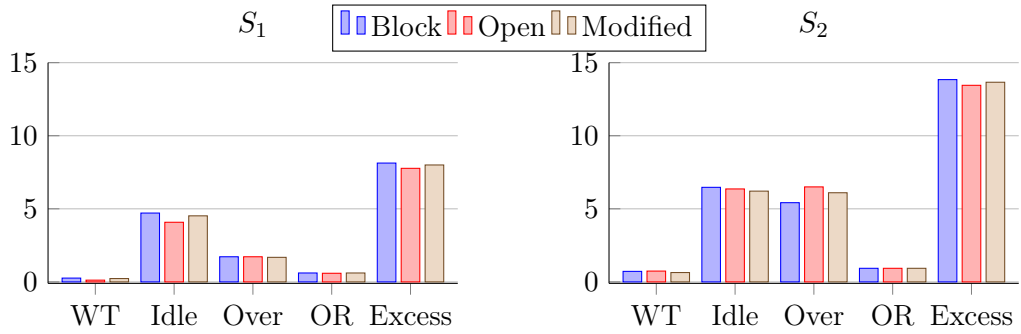

Looking at Table \ref{tab:unc:policy_00} and Figure \ref{fig:uncertainty_policy_00}, we can say that if we set $\Theta_1=\Theta_2=0$, the three policies exhibit very similar performance, indicating that the system behavior is largely insensitive to the scheduling policy under these conditions. This suggests that, under low uncertainty conditions, the system is primarily constrained by shared structural resources rather than by the specific scheduling rule.
The number of excess beds required in $S_2$ is higher than in $S_1$, likely because more patients are scheduled under a lower LoS, which leads to higher variability. Also the increasing of idle time and overtime is explained by the difference in the number of scheduled patients.

\begin{table}[tbp]
    \small
    \centering
    \caption{Policy Comparison with $\Theta_1=50$ and $\Theta_2=65/90$ under uncertainty for instance $I_2$. The results are computed only for the instances that are feasible for each policy.}
    \label{tab:unc:policy_5095} 
    \renewcommand{\arraystretch}{1.2}
    \resizebox{\textwidth}{!}{%
    \begin{tabular}{l r r r r r r r r r r r r r r}
        & &  \multicolumn{6}{c}{$S_1$} & & \multicolumn{6}{c}{$S_2$} \\
        \cmidrule(lr){3-8} \cmidrule(lr){10-15}
        \textbf{Policy} & &
        \textbf{WT} &
        \textbf{Idle} & \textbf{Over} & \textbf{OR real} & \textbf{Excess} & \textbf{Beds (\%)}& &
        \textbf{WT}  &
        \textbf{Idle} & \textbf{Over}& \textbf{OR real} & \textbf{Excess} & \textbf{Beds (\%)}\\
         \midrule
        Block && 0.49 & 6.55 & 3.39 & 0.73 & 5.66 & 0.95 && 0.87 & 6.67 & 7.40 &0.99& 8.84 & 0.90 \\
Open  && 0.59 & 7.34 & 4.88 & 0.72&5.62 & 0.95 &  & 0.97 & 7.51 & 8.16 & 1.00 &8.42 & 0.87 \\
Modified && 0.52 & 6.75 & 4.03 &0.73& 5.43 & 0.95 & & 0.88& 7.04 & 7.44 & 1.00& 8.59 & 0.90 \\
        \bottomrule
    \end{tabular}}
\end{table}

\begin{figure}[h!]
\caption{Comparison of uncertainty indicators across the three scheduling policies, considering $\Theta_1 = 50$ and $\Theta_2=65/90$ for instance $I_2$.}
\label{fig:uncertainty_policy5090}
\centering
\begin{tikzpicture}

\begin{groupplot}[
group style={
    group size=2 by 1,
    horizontal sep=2cm,
},
enlarge x limits=0.15,
ybar,
width=7cm,
height=4.5cm,
symbolic x coords={WT,Idle,Over,OR,Excess},
xtick=data,
ymin=-0.1,
ymax=10,
ymajorgrids=true,
axis x line*=bottom,
axis y line*=left,
legend style={
    at={(1.15,1.05)},
    anchor=south,
    legend columns=3,
},
]

\nextgroupplot[
title={$S_1$},
ybar,
bar width=7pt,
]

\addplot coordinates {
(WT,0.49)
(Idle,6.55)
(Over,3.39)
(OR,0.73)
(Excess,5.66)
};

\addplot coordinates {
(WT,0.59)
(Idle,7.34)
(Over,4.88)
(OR,0.72)
(Excess,5.62)
};

\addplot coordinates {
(WT,0.52)
(Idle,6.75)
(Over,4.03)
(OR,0.73)
(Excess,5.43)
};

\legend{Block,Open,Modified}

\nextgroupplot[
title={$S_2$},
ybar,
bar width=7pt,
]

\addplot coordinates {
(WT,0.87)
(Idle,6.67)
(Over,7.40)
(OR, 0.99)
(Excess,8.84)
};

\addplot coordinates {
(WT,0.97)
(Idle,7.51)
(Over,8.16)
(OR, 1.00)
(Excess,7.42)
};

\addplot coordinates {
(WT,0.88)
(Idle,7.04)
(Over,7.44)
(OR, 1.00)
(Excess,8.59)
};

\end{groupplot}

\end{tikzpicture}
\end{figure}
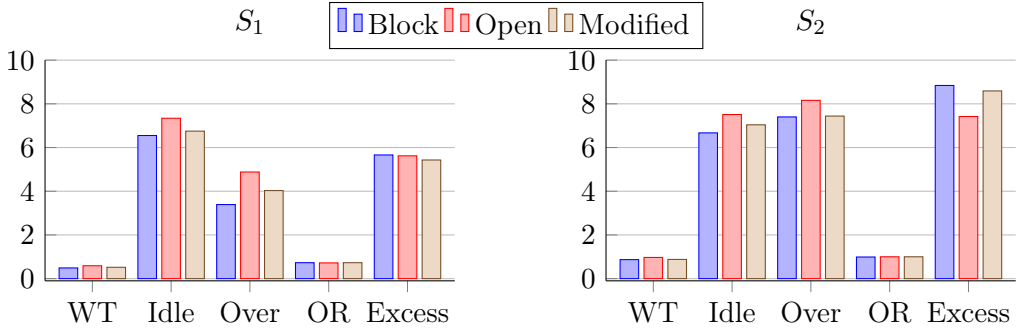

From Table~\ref{tab:unc:policy_5095} and Figure~\ref{fig:uncertainty_policy5090}, we observe that, also in this case, the differences between policies are small. Although the open scheduling policy performs worst in terms of objective function, it also performs worst in terms of uncertainty-related measures, generally yielding the highest waiting time, idle time, and overtime. This suggests that, while this policy is more flexible and may improve patient accommodation, it can also lead to increased delays and overtime.
\section{Conclusions}

This paper addressed the integrated operating room scheduling problem under weekend bed capacity constraints. We proposed a column generation framework that jointly solves the MSSP and the SCAP, while simultaneously determining the assignment of specialties, surgeons, and patients to OR blocks, sequencing surgeries, and accounting for downstream bed occupancy. The resulting optimization model maximizes the total priority of scheduled patients while controlling weekend bed occupancy and promoting efficient OR utilization.

An extensive computational study was conducted under different operating settings, scheduling policies, and resource availability scenarios. The proposed CG approach consistently outperformed the compact MILP formulation for the majority of instances. Although the MILP approach remained competitive for some small instances, the CG approach obtained better solutions within shorter computational times for the remaining cases. For larger instances, the MILP approach was unable to find feasible solutions within the imposed time limit. In contrast, the proposed approach consistently generated high-quality feasible solutions for all tested large instances.

The computational analysis highlighted several managerial insights. First, weekend bed availability has a significant impact on the number of scheduled patients, showing that downstream bed capacity may become the primary bottleneck of the surgical scheduling process. Increasing weekend bed availability improves patient access but requires additional staffing resources. Second, imposing a minimum OR utilization requirement increases OR efficiency but may reduce the overall priority of the scheduled patients by forcing the inclusion of lower-priority surgeries. Finally, among the scheduling policies considered, modified block scheduling provides the best compromise between scheduling flexibility and computational efficiency, achieving the best overall performance in most of the tested instances.

The robustness analysis under uncertainty further demonstrated that the proposed schedules remain stable when surgery durations and patient LoS deviate from their expected values. Waiting times, idle times, and overtime remain limited, while OR and bed utilization exhibit only minor variations, confirming the practical applicability of the proposed approach.

Future research may focus on strengthening the column generation framework through more advanced pricing strategies and branching rules, as well as extending the model to explicitly incorporate uncertainty within the optimization process and additional operational constraints arising in real hospital environments.

\section*{Acknowledgments}
The work of the author Sara Cambiaghi is supported by the Ph.D.\ fellowship funded under the Italian National Recovery and Resilience Plan (PNRR), as part of the program for research doctorates and innovative doctorates for public administration and cultural heritage (D.M.\ 118, 02-03-2023).
\appendix
\section{Proof of Lemma 1}
\label{Appendix2}

First suppose that $|I| = 3$. If all patients are assigned to the same surgeon or to three different surgeons, there is nothing to prove. Thus, assume that two patients, $i$ and $i'$, are assigned to surgeon $s_1$, and the remaining patient $j$ is assigned to surgeon $s_2$.  

Suppose, by contradiction, that the sequence $(iji')$ dominates $(ii'j)$. This would require:
\begin{equation}
    \begin{cases}
        \tilde{\pi_i} + \tilde{\pi_j} + \tilde{\pi_{i'}} - d^2_{s_1} - d^2_{s_2} - d^3 - d^4 
        \ge \tilde{\pi_i} + \tilde{\pi_{i'}} + \tilde{\pi_j} - d^2_{s_1} - d^2_{s_2} - d^3 - d^4, \\[3pt]
        \delta_i + \tau_{ij} + \delta_j + \tau_{ji'} + \delta_{i'}
        \le \delta_i + \tau_{ii'} + \delta_{i'} + \tau_{i'j} + \delta_{j},
    \end{cases}
\end{equation}
with at least one of the two inequalities being strict.  
The first inequality is actually an equality, so $(iji')$ dominates $(ii'j)$ only if
\[
    \delta_i + \tau_{ij} + \delta_j + \tau_{ji'} + \delta_{i'}
    <
    \delta_i + \tau_{ii'} + \delta_{i'} + \tau_{i'j} + \delta_{j}
    \quad\Longleftrightarrow\quad
    \tau_{ij} + \tau_{ji'} < \tau_{ii'} + \tau_{i'j}.
\]
By assumption~\textit{a}, $\tau_{ii'} \le \tau_{ij}$, so the inequality reduces to
\(
    \tau_{ji'} < \tau_{i'j},
\)
which is impossible by assumption~\textit{b}.

Now consider the case $|I| = 4$. We distinguish several assignments:

\begin{itemize}

    \item If all patients are assigned to the same surgeon or to four different surgeons, there is nothing to prove.
    \item Patients $i,i',i''$ assigned to $s_1$ and $j$ assigned to $s_2$:
    \begin{itemize}
        \item $i\underbracket{i' j i''}_{>\, i' i'' j} > i i' i'' j$;
        \item $\underbracket{i j i'}_{>\, i i' j} i'' > i i' j i'' > i i' i'' j$.
    \end{itemize}
    \item Patients $i,i'$ assigned to $s_1$, and $j,j'$ assigned to $s_2$:
    \begin{itemize}
        \item $j\underbracket{i j' i'}_{>\, i i' j} > j i i' j$.  
        Define $\hat{i} = ii'$ as a single “combined’’ patient with surgery duration  $\delta_i+\tau_{ii'}+\delta_{i'}$. Then $j i i' j' = j \hat{i} j' > j j' \hat{i} = j j' i i'.$
        \item Same argument applies to the case $jii'j'$.
    \end{itemize}
    \item Patients $i,i'$ assigned to $s_1$, $j$ to $s_2$, and $z$ to $s_3$:
    \begin{itemize}
        \item $z\underbracket{i j i'}_{>\, i i' j} > z i i' j$;
        \item $i z j i' = i \hat{zj} i' > i i' \hat{zj} = i i' z j$,  
        where $\hat{zj}$ is treated as a single patient with surgery duration  
        $\delta_z + \tau_{zj} + \delta_j$.
    \end{itemize}
\end{itemize}

By induction, the claim holds for any cardinality of the set $I$.


\singlespacing
\bibliographystyle{elsarticle-harv}
\bibliography{ref_short}

\section*{Supplementary Material}
This supplementary material provides additional mathematical formulations, as well as examples of instances and corresponding solutions. In particular, Section \ref{Appendix4} presents the mathematical formulation of the subproblem, Section \ref{Appendix5} provides an example of an instance, and Section \ref{Appendix1} presents an example of a solution considering block, open, and modified scheduling policies.
 
\section{MILP formulation of the subproblem}
\label{Appendix4}
\label{section:cg:sp:milp}
In the following, we present the mathematical formulation of the subproblem. Given an OR block $(r,e)$ and, forthe dedicated OR blocks, a selected specialty $j$, our goal is to determine the optimal schedule for OR block $(r,e)$ (and, for the dedicated OR blocks, the patients belonging to specialty $j$).

The decision variables are defined as follows:
\begin{itemize}
\item $\alpha_i = 1$ if patient $i$ is included in the schedule, and $0$ otherwise,
\item $o_{ii'} = 1$ if patient $i$ is scheduled immediately before patient $i'$, and $0$ otherwise,
\item $f_i = 1$ if patient $i$ is the first patient in the schedule, and $0$ otherwise,
\item $l_i = 1$ if patient $i$ is the last patient in the schedule, and $0$ otherwise,
\item $t_i \ge 0$ is the start time of the surgery of patient $i$, and is set to $0$ if the patient is not scheduled.
\item $z_j = 1$ if specialty $j$ is selected in the schedule,
\item $\rho_s = 1$ if surgeon $s$ is selected in the schedule, and $0$ otherwise,
\item $b_{id'} = 1$ if patient $i$ occupies a bed on day $d'$, and $0$ otherwise,
\item $\psi_{jd'} \ge 0$ is the number of beds occupied on day $d'$ for specialty $j$.
\end{itemize}

We now present all the constraints of the model. For notational simplicity, the constraints are defined for all $r \in \mathcal{R}$ and $e \in \mathcal{E}$, but are enforced only for pairs $(r,e) \in \mathcal{B}$.

 If a patient is not assigned to the OR block $(r, e)$, the corresponding start-time variable is set to zero.
\begin{equation}
    t_i \le M\alpha_i \quad \forall i\in \mathcal{I}
\end{equation}

Only the selected specialty is eligible:
\begin{equation}
\begin{cases}
    \sum_{j'\in \mathcal{J}} z_j' = 1\\
    z_{j} = 1
\end{cases}
\end{equation}
This constraint is enforced only for dedicated OR blocks. For these blocks, the problem structure allows the subproblem to be decomposed by specialty. Specifically, for each OR block, the subproblem is solved separately for each specialty compatible with that block. This constraint ensures that only the selected specialty is considered in the corresponding subproblem.

A patient can be scheduled in the OR block $(r,e)$ only if that block is assigned to the patient’s specialty $j_i$.
\begin{equation}
   \sum_{i\in \mathcal{I}\setminus \mathcal{I}^j} \alpha_i = 0
\end{equation}

A patient can be scheduled in the OR block $(r,e)$ only if their assigned surgeon $s_i$ is also assigned to that OR block.
Each surgeon can be assigned to the OR block  $(r,e)$ only if at least one patient assigned to them is scheduled in that block.
\begin{align}
        \alpha_i & \le \rho_{s_i} \quad \forall i\in \mathcal{I} \\
        \rho_s & \le \sum_{i\in \mathcal{I}^s}\alpha_i \quad \forall s\in \mathcal{S}
\end{align}

Each surgery must be completed before the end of the time available for the OR block.
\begin{equation}
   t_i + \delta_i \le \xi_{re} + M(1-\alpha_i) \quad \forall i\in \mathcal{I}
\end{equation}

This constraint defines the start times of surgeries and ensures that procedures scheduled in the
OR block do not overlap.
\begin{equation}
    t_{i'} \ge t_i + \delta_i + \tau_{ii'} -M(1-o_{ii'}) \quad \forall i\neq i'\in \mathcal{I}
\end{equation}

Exactly one patient is designated as the first and one as the last in the schedule. Since the subproblem is solved for each OR block, we enforce that there is exactly one first patient and exactly one last patient. Consequently, a feasible schedule is generated for each block, and the master problem determines how many schedules to select.
\begin{align}
    \sum_{i\in \mathcal{I}} f_i &= 1  \\
    \sum_{i\in \mathcal{I}} l_i &= 1 
\end{align}

A patient assigned as the first, last, or intermediate patient in the OR block sequence must be scheduled in that block.
\begin{align}
        f_i & \le \alpha_i & \quad \forall i\in \mathcal{I}\\
        l_i & \le \alpha_i & \quad \forall i\in \mathcal{I}\\
        o_{ii'} & \le \alpha_i & \quad \forall i\in \mathcal{I} \\
        o_{i'i} &\le \alpha_i & \quad \forall i\in \mathcal{I}
\end{align}

A scheduled patient is either the first patient in the sequence or has a predecessor.
\begin{equation}
    \alpha_i = f_i + \sum_{i'\in \mathcal{I}, i'\neq i}o_{i'i} \quad \forall i\in \mathcal{I} 
\end{equation}

A scheduled patient is either the last patient in the sequence or has a successor.
\begin{equation}
    \alpha_i = l_i + \sum_{i'\in \mathcal{I}, i'\neq i}o_{ii'} \quad \forall i\in \mathcal{I}
\end{equation}

\begin{equation}
    t_i \le (1-f_i)M_e \quad \forall i\in \mathcal{I} 
\end{equation}

Let $d$ be the unique day such that $e\in \mathcal{E}^d$. We define $b_{id'}=1$ if patient $i$ occupies a bed on day $d'$. Then:
\begin{equation}
    \sum_{d'=d}^{\min\{d+\nu_i, |\mathcal{D}|\}} b_{id'} \ge \min\{\nu_i, |\mathcal{D}| - d\} \alpha_i \quad \forall i\in\mathcal{I}
\end{equation}

Next, we define the bed occupancy for each specialty $j'$ on day $d'$ as:
\begin{equation}
    \psi_{j'd'} = \sum_{i\in \mathcal{I}|\gamma_{ij'}=1} b_{id'} \quad \forall j'\in \mathcal{J}
\end{equation}
Let us note that, for the OR blocks under the block scheduling policy, it is sufficient to define this constraint only for the considered specialty $j$; for all other specialties $j' \in \mathtt{J} \setminus \{j\}$, we have $\psi_{j'd} = 0$.


Our objective is to maximize the following expression:
\begin{align}
    \max \quad & 
    \sum_{i\in \mathcal{I}}\pi_i\alpha_{i} -
    \sum_{i\in \mathcal{I}} \alpha_{i} \mu^1_i - 
    \sum_{s\in \mathcal{S}}\rho_{s}\mu^2_{se} - 
    \mu^3_{re} - \\
    &\quad \sum_{j\in \mathcal{J}}\sum_{d'\in\mathcal{D}}\psi_{jd'}\mu^4_{jd'} - 
    \sum_{j\in \mathcal{J}} \sum_{d'\in\mathcal{D}^w} \psi_{jd'} \mu^5_{jd'} 
    - \left\{\sum_{i\in \mathcal{I}}\left[\alpha_i(\delta_i+\tilde{\tau})\right]-\tilde{\tau}\right\} \mu^6
\end{align}
\section{Example of instances}
\label{Appendix5}
\begin{figure}[H]
   \centering
    \caption{Example of instance $I_1$, setting $S_1$. From the first plot on the left, we show the following components: \textbf{Specialty--Patients}, representing the number of patients in the waiting list for each specialty; \textbf{Surgery Duration per Specialty}, representing the distribution of surgery durations for each specialty; and \textbf{LoS per Specialty}, representing the LoS distribution for each specialty (in setting $S_2$, this is a discrete uniform distribution between 0 and 4). \textbf{Beds per Specialty} represents the number of beds assigned to each specialty. \textbf{Turnover Time Distributions} represent the distribution of turnover times, distinguishing between turnovers for patients assigned to the same surgeon and those assigned to different surgeons.  \textbf{Surgeon--Patients} represents the number of patients in the waiting list for each surgeon.}
       \label{fig:instance:variable}
   \includegraphics[width=\textwidth]{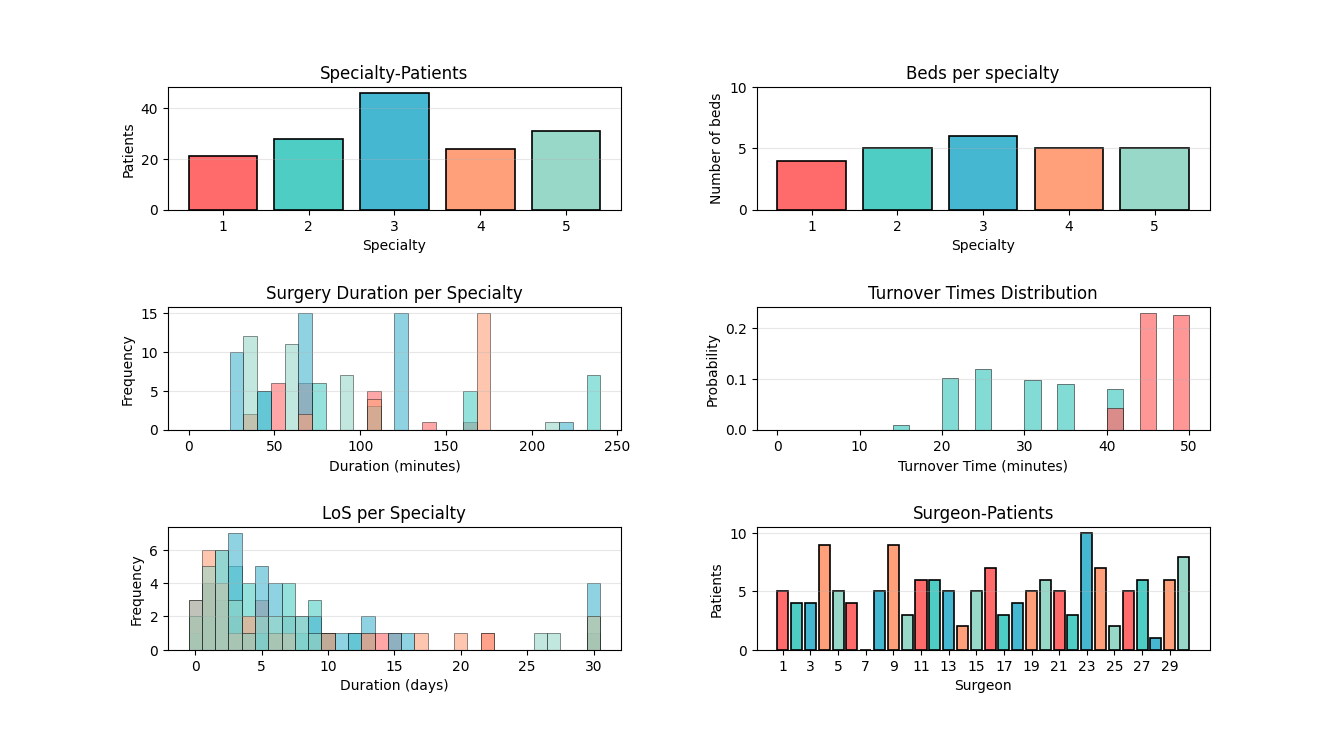}
      
   \end{figure}

In the following, an overview of an example instance is provided in Figure \ref{fig:instance:variable}. 
\section{Examples of solutions}
\label{Appendix1}

\begin{figure}[H]
\caption{Example of solutions with the block scheduling policy}
    \label{fig:example:block}
    \centering
    \begin{subfigure}{0.8\textwidth}
        \centering
        \includegraphics[width=\linewidth]{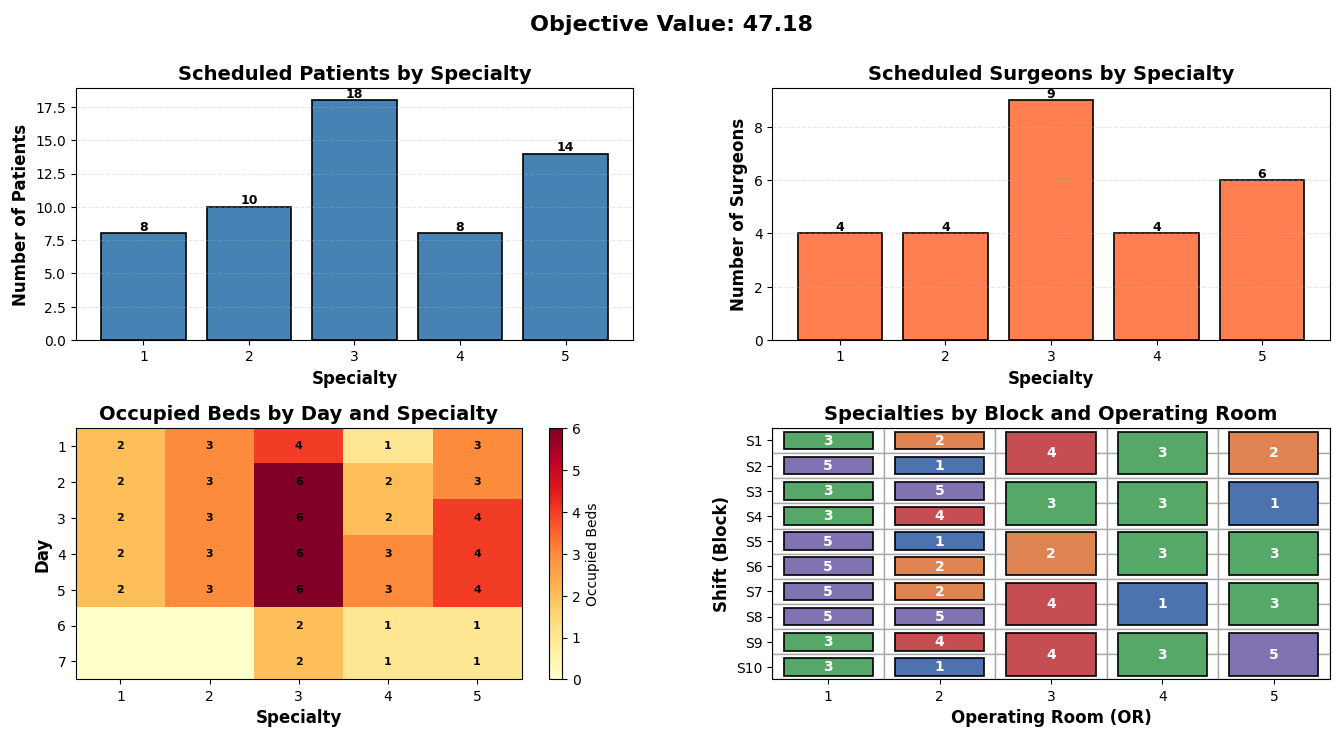}
        \caption{Setting $S_1$}
    \end{subfigure}
    \vspace{0.5cm}
    \begin{subfigure}{0.8\textwidth}
        \centering
        \includegraphics[width=\linewidth]{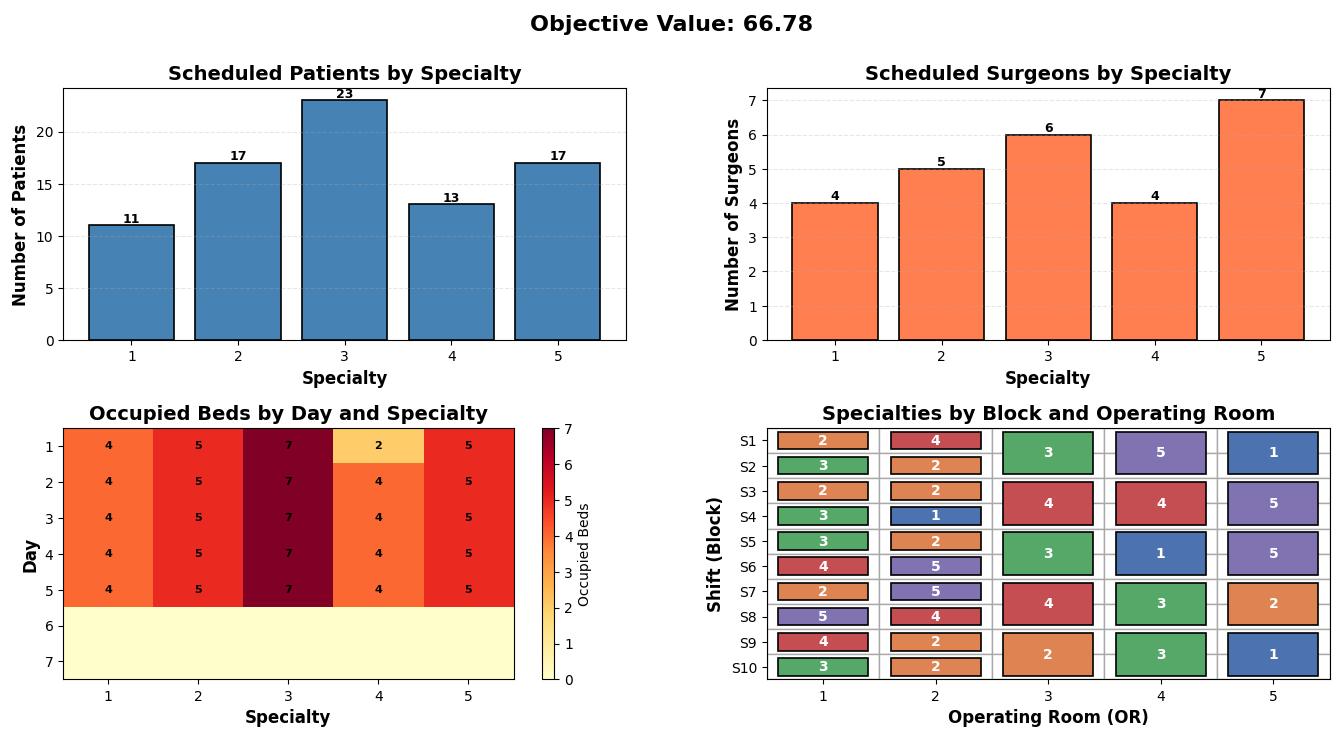}
        \caption{Setting $S_2$}
    \end{subfigure}
    
\end{figure}

\begin{figure}[H]
\caption{Example of solutions with the open scheduling policy}
    \label{fig:example:open}
    \centering
    \begin{subfigure}{0.8\textwidth}
        \centering
        \includegraphics[width=\linewidth]{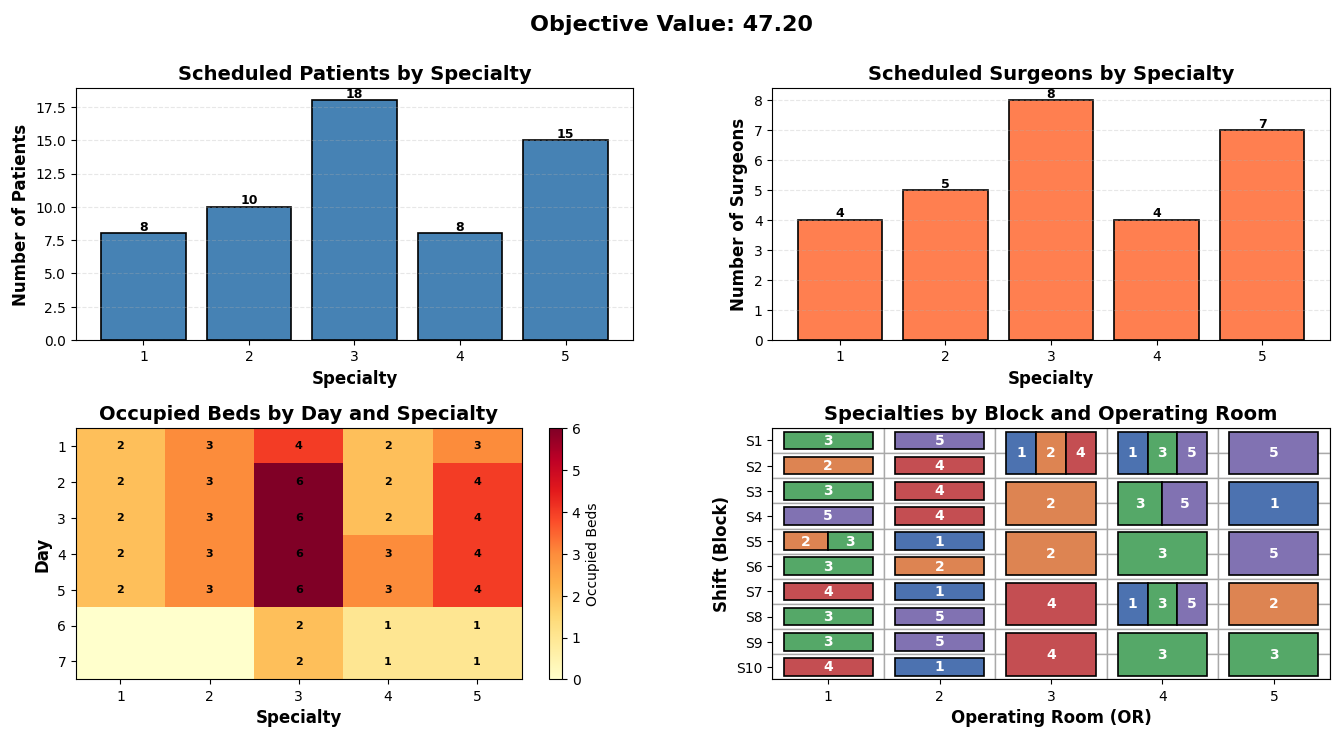}
        \caption{Setting $S_1$}
    \end{subfigure}
    \vspace{0.5cm}
    \begin{subfigure}{0.8\textwidth}
        \centering
        \includegraphics[width=\linewidth]{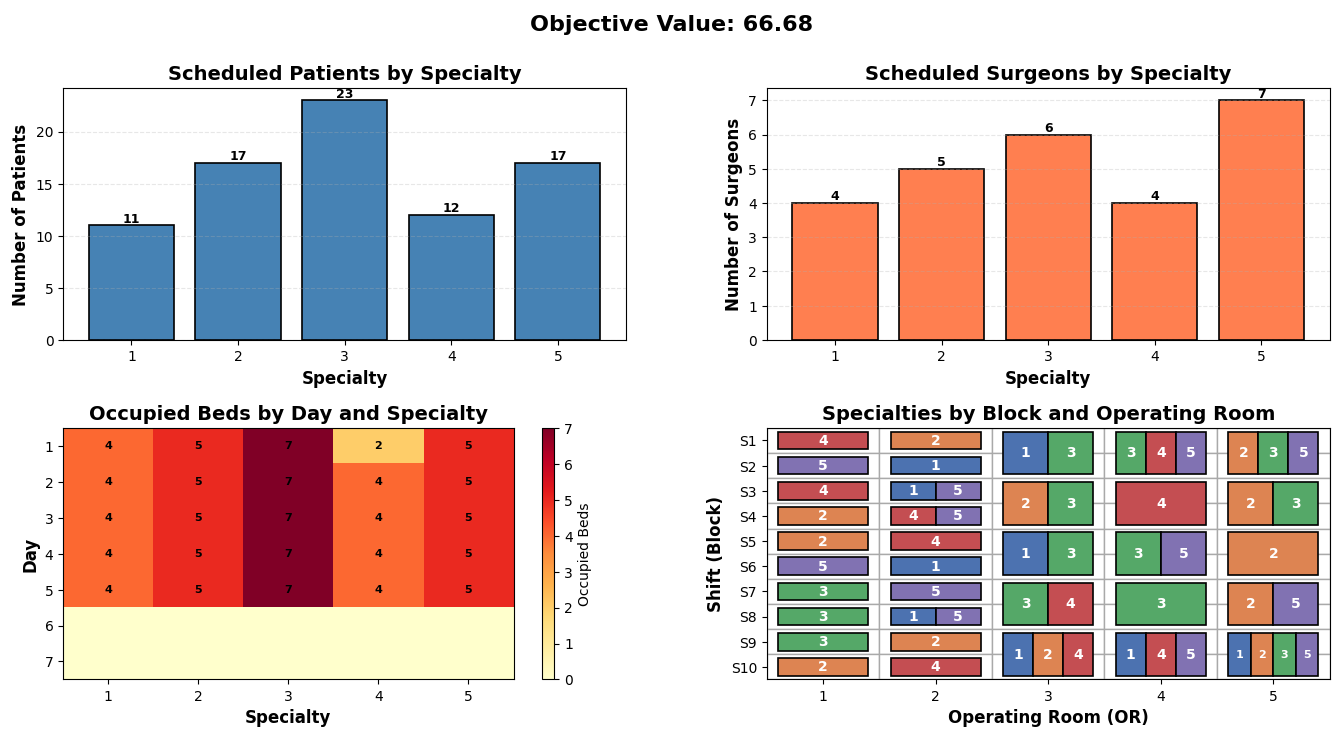}
        \caption{Setting $S_2$}
    \end{subfigure}
\end{figure}

\begin{figure}[H]
\caption{Example of solutions with the modified scheduling policy}
    \label{fig:example:modified}
    \centering
    \begin{subfigure}{0.8\textwidth}
        \centering
        \includegraphics[width=\linewidth]{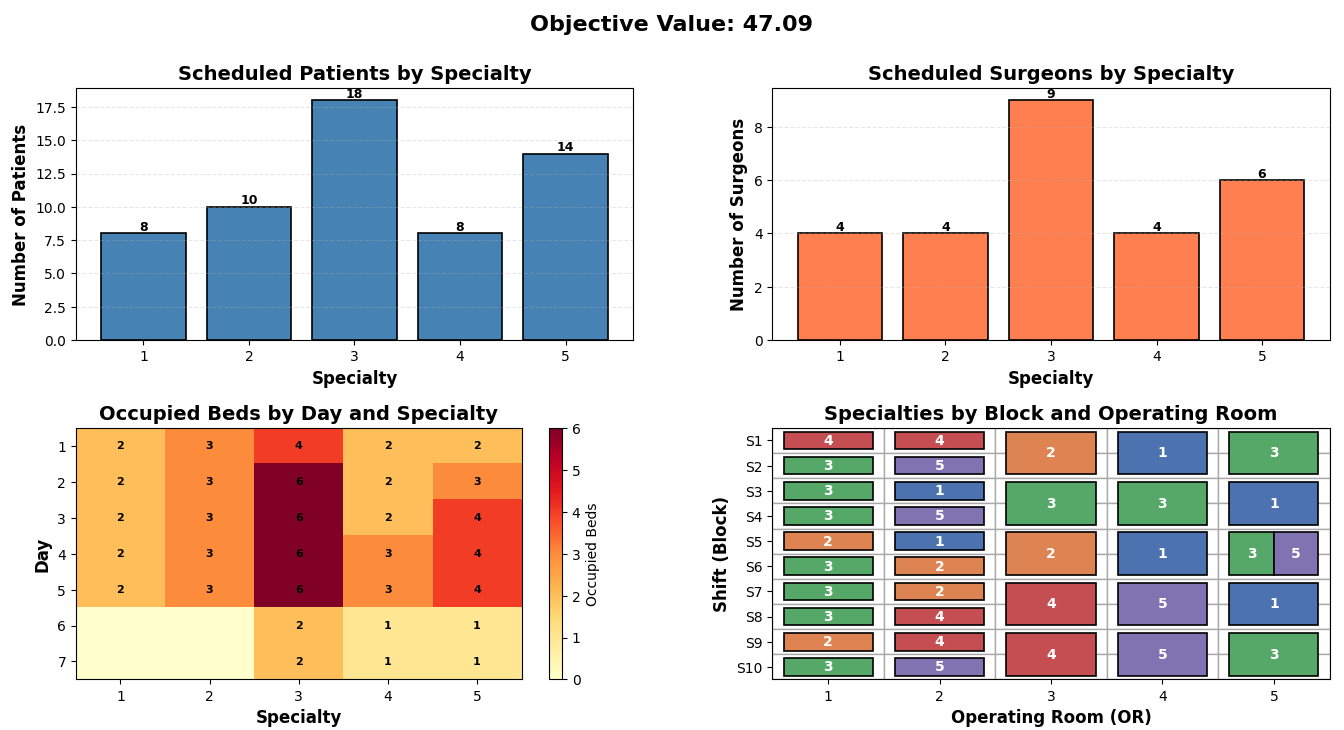}
        \caption{Setting $S_1$}
    \end{subfigure}
    \vspace{0.5cm}
    \begin{subfigure}{0.8\textwidth}
        \centering
        \includegraphics[width=\linewidth]{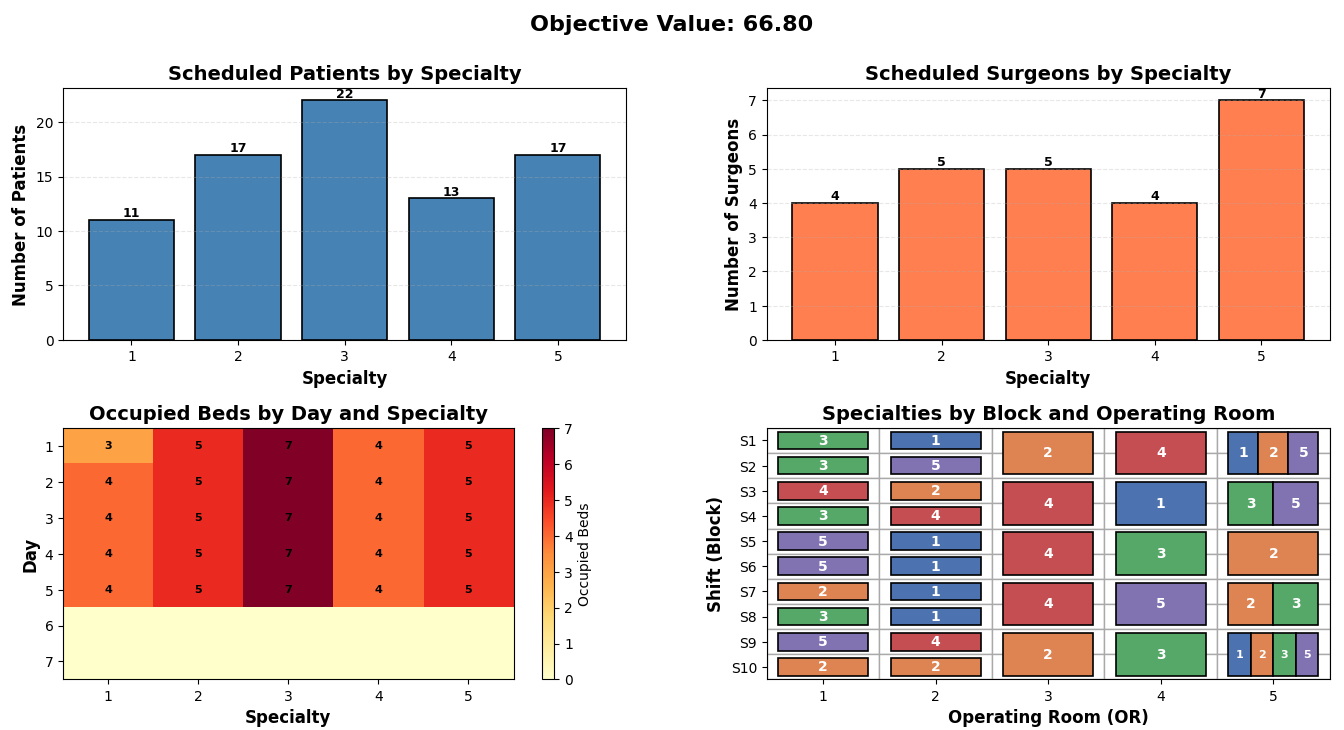}
        \caption{Setting $S_2$}
    \end{subfigure}
\end{figure}

In the following, we report examples of solutions for the instance $I_1$ under the two settings $S_1$ and $S_2$, using the block scheduling (Figure \ref{fig:example:block}), the open scheduling (Figure \ref{fig:example:open}) and the modified scheduling (Figure \ref{fig:example:modified}) policies.

\end{document}